\documentclass[12pt, reqno]{amsart}
\usepackage{mathrsfs}
\usepackage[active]{srcltx}
\usepackage{mathrsfs}
\usepackage{color}

\usepackage{color}
\usepackage[unicode]{hyperref}
\hypersetup{
    colorlinks = true,%
    citecolor = [rgb]{0.0,0.0,0.9},
    filecolor=black,%
    linkcolor = [rgb]{0.65,0.0,0.0},%
    anchorcolor = red,
    pagecolor = red,
    urlcolor= [rgb]{0.65,0.0,0.0}
}

\newtheorem{example}{ Example}[section]
\newtheorem{proposition}{Proposition}[section]
\newtheorem{theorem}{Theorem}[section]
\newtheorem{lemma}{Lemma}[section]
\newtheorem{corollary}{Corollary}[section]

\newtheorem{remark}{Remark}[section]

\begin{document}

\title[Moser-Trudinger inequalities on Riemannian manifolds]
{New geometric aspects of Moser-Trudinger inequalities on Riemannian manifolds: the non-compact case}

\author{Alexandru Krist\'aly}
\address{Department of Economics, Babe\c s-Bolyai University, Str. T. Mihali 58-60, 400591
	Cluj-Napoca, Romania \&  Institute of Applied Mathematics, \'Obuda University,
	B\'ecsi \'ut 96, 1034 Budapest, Hungary}
\email{alex.kristaly@econ.ubbcluj.ro; kristaly.alexandru@nik.uni-obuda.hu}
\thanks{Research supported by the National Research, Development and Innovation Fund of Hungary, financed under the K$\_$18 funding scheme, Project No.  127926, and by the STAR-UBB Institute, Cluj-Napoca, Romania.}





\begin{abstract} In the first part of the paper we investigate some geometric features of Moser-Trudinger inequalities on
complete non-compact  Riemannian manifolds. By exploring
rearrangement arguments, isoperimetric estimates, and gluing local
uniform estimates via Gromov's covering lemma, we provide a Coulhon, Saloff-Coste and Varopoulos type characterization concerning the
validity of Moser-Trudinger inequalities on complete non-compact
$n-$dimensional Riemannian manifolds $(n\geq 2)$ with Ricci
curvature bounded from below. Some sharp consequences are also presented both for non-negatively and non-positively curved Riemannian manifolds, respectively. In the second part, by combining variational arguments and a Lions type symmetrization-compactness principle, 
we guarantee the existence of a non-zero isometry-invariant solution
for an elliptic problem involving the $n-$Laplace-Beltrami operator and a critical nonlinearity on $n-$dimensional homogeneous 
 Hadamard manifolds. Our results comp\-lement in several directions those of Y. Yang [J. Funct. Anal., 2012]. 
\end{abstract}

\maketitle

\begin{center}
\textit{Dedicated to my children, Mar\'ot, Bora, Zonga and Bendeg\'uz.} 
\end{center}


 \section{\sc Introduction}
\subsection{Objectives} The Moser-Trudinger inequality, as the borderline case of
Sobolev inequalities, plays a crucial role in the theory of
geometric functional analysis and its applications in the study of
quasilinear elliptic problems on the Sobolev space $W^{1,n}$ defined
on $n-$dimensional geometric objects,  $n\geq 2$.

In the present paper we investigate the influence of geometry of
complete non-compact Riemannian manifolds to the validity, sharpness
and further aspects of Moser-Tru\-dinger inequalities. Roughly
speaking, we shall
\begin{itemize}
  \item characterize the validity of  Moser-Trudinger inequalities on complete non-compact Riemannian manifolds with  Ricci curvature bounded from
  below in terms of the volume growth of geodesic balls  (no assumption  on the injectivity radius is required);
  \item provide sharp consequences both on non-negatively and non-positively curved Riemannian manifolds;
  \item guarantee the existence of a non-zero isometry-invariant solution for a quasilinear elliptic problem  on $n-$di\-men\-sional homogeneous Hadamard
manifolds which involves the $n-$Laplace-Beltrami operator and a
term  with critical  growth.
\end{itemize}
Before to state our results we recall some features of
the Moser-Trudinger inequality which will be used in the sequel.

\subsection{Short overview of Moser-Trudinger
inequalities}\label{sect1.1} Let $\Omega$ be an open subset of the
Euclidean space $ \mathbb R^n$ $(n\geq 2)$ with finite Lebesgue
measure. It is well known that the borderline case of the Sobolev
embeddings $W_0^{1,p}(\Omega)\hookrightarrow L^q(\Omega)$, where
$1\leq q\leq \frac{np}{n-p}$ and $1<p<n$, has an unusual
behavior; indeed, when $n=p$, the Sobolev space $W_0^{1,n}(\Omega)$
cannot be continuously embedded  into $L^\infty(\Omega)$, although
formally this should be the case. Motivated by this phenomenon,
Trudinger \cite{Trudinger} proved that
$W_0^{1,n}(\Omega)\hookrightarrow L_{\psi_n}(\Omega),$ where
$L_{\psi_n}(\Omega)$ is the Orlicz space associated with the Young
function $\psi_n(s)=e^{\alpha |s|^\frac{n}{n-1}}-1$ for  $\alpha>0$
sufficiently small. A few years later,
 Moser \cite{Moser}
stated the sharp version of this embedding, by proving that there
exists $M_0=M_0(n)>0$ depending only on $n$ such that
\begin{equation}\label{Moser-sharp}
   \sup_{u\in \mathcal H}\int_\Omega e^{\alpha|u|^\frac{n}{n-1}}{\rm d}x=\left\{ \begin{array}{lll}
 M_0 {\rm Vol}_e(\Omega) &\mbox{if} &  \alpha\in [0,\alpha_n]; \\
+\infty &\mbox{if} &  \alpha>\alpha_n;
 \end{array}\right.
\end{equation}
here $\mathcal H=\left\{u\in W_0^{1,n}(\Omega): 
\int_\Omega |\nabla u|^n{\rm d}x\leq 1\right\},$ Vol$_e(\cdot)$ is
the Euclidean volume, $\omega_{n-1}$ is the area of the unit sphere
$\mathbb S^{n-1}\subset\mathbb R^n$ and
$$\alpha_n=n\omega_{n-1}^\frac{1}{n-1}$$ is the {\it critical exponent}.

The Moser-Trudinger inequality (\ref{Moser-sharp}) became in this
way the starting point of further studies in various directions,
both in the Euclidean and non-Euclidean settings. In the Euclidean
case, milestone results can be found concerning the sharpness and
existence of extremal functions for the classical Moser-Trudinger
inequality both on bounded and unbounded sets, see e.g. Carleson and
Chang \cite{CC}, Flucher \cite{Flucher}, Lin \cite{Lin}, Li and Ruf
\cite{LR}. In particular, if $n\geq 2$ and
\begin{equation}\label{phi-def}
    \Phi_n(t)=e^t-\sum_{k=0}^{n-2}\frac{t^k}{k!},
\end{equation}
Li and Ruf \cite{LR} proved that
\begin{equation}\label{Li-Ruf-eredmeny}
    S_n^{LR}:=\sup_{u\in W^{1,n}(\mathbb R^n),\ \|u\|_{0,1}\leq 1}\int_{\mathbb
    R^n}\Phi_n(\alpha_n|u|^\frac{n}{n-1}){\rm d}x<\infty,
\end{equation}
where $\|u\|_{0,1}^n=\int_{\mathbb R^n}(|\nabla
u|^n+|u|^n){\rm d}x.$ 
The constant $\alpha_n$ in (\ref{Li-Ruf-eredmeny}) is sharp;
although the integral in (\ref{Li-Ruf-eredmeny}) is finite for every
$\alpha>0$ instead of $\alpha_n$, the supremum is infinite for
$\alpha>\alpha_n.$ Improvements and higher order extensions of the
Moser-Trudinger inequality can be found e.g. in Adams \cite{Adams},
Adimurthi and Druet \cite{Adi-Druet}, Cianchi, Lutwak, Yang and
Zhang \cite{Cianchi}, Ibrahim, Masmoudi and Nakanishi \cite{IMN},
Masmoudi and Sani \cite{Mas-Sani}, Ruf and Sani
\cite{Ruf-Sani-TAMS}, and references therein.

Moser-Trudinger inequalities in the non-Euclidean setting captured
also special attention. On one hand, sharp Moser-Trudinger
inequalities are established in Heisenberg and Carnot groups, see
Cohn and Lu \cite{CL}, Lam and Lu \cite{LamLu-Adv}, Balogh, Manfredi
and Tyson \cite{BMT}, and on CR spheres, see Branson, Fontana and
Morpurgo \cite{BFM}. On the other hand, deep achievements can be
found in the study of Moser-Trudinger inequalities on Riemannian
manifolds which are particularly important from the viewpoint  of
the present paper.

   Let $n\geq 2$  and $(M,g)$ be an
$n-$dimensional Riemannian manifold endowed with its canonical
 volume form ${\rm d}v_g$. For $\tau>0$ fixed, on the usual
Sobolev space $W^{1,n}(M)=W_0^{1,n}(M)$ defined on $(M,g)$, see
Hebey \cite{Hebey}, we consider the equivalent norms
$$\|u\|_{0,\tau}=\left(\|\nabla_g u\|^n_{L^n(M)}+\tau^n \|u\|^n_{L^n(M)}\right)^\frac{1}{n}\ {\rm and}\ \|u\|_{1,\tau}=\|\nabla_g u\|_{L^n(M)}+\tau\|u\|_{L^n(M)},$$
where the Lebesgue norms $\|\cdot\|_{L^n(M)}$ are defined by means
of the volume form ${\rm d}v_g$.  According to these norms, for
every $\alpha> 0$, $\tau>0$ and $i\in \{0,1\}$, we introduce the
quantities
 $$   S_{\alpha,\tau}^i(M,g):=\sup_{u\in W^{1,n}(M),\ \|u\|_{i,\tau}\leq
1}\int_M \Phi_n(\alpha|u|^\frac{n}{n-1}){\rm d}v_g,$$
where $\Phi_n$  is from (\ref{phi-def}).  Since $\|u\|_{0,\tau}\leq
\|u\|_{1,\tau}$, then $S_{\alpha,\tau}^1(M,g)\leq
S_{\alpha,\tau}^0(M,g)$ for every $\alpha,\tau> 0.$

Let $i\in \{0,1\}.$ If there exists $\alpha> 0$ and $\tau>0$ such
that $S_{\alpha,\tau}^i(M,g)<+\infty$, we say that the {\it
Moser-Trudinger inequality $({\bf MT})_{\alpha,\tau}^i$ holds on
$(M,g)$}. Contrary, if $S_{\alpha,\tau}^i(M,g)=+\infty$ for some
$\alpha> 0$ and $\tau>0$, we say that the {\it Moser-Trudinger
inequality $({\bf MT})_{\alpha,\tau}^i$ fails on $(M,g)$}.

On one hand, when $(M,g)$ is an $n-$dimensional {\it compact}
Riemannian manifold {\it without boundary}, then for every
$\alpha\in [0,\alpha_n]$ and $\tau>0$, the Moser-Trudinger
inequality $({\bf MT})_{\alpha,\tau}^0$ holds on $(M,g)$ and the
critical exponent $\alpha_n$ is sharp,  see Li \cite{Li}; a higher
order extension of Li's result can be found in do \'O and Yang
\cite{doO-Yang}. Note that both papers \cite{doO-Yang} and \cite{Li}
are extensions of Fontana \cite{Fontana} replacing the constraints
$\int_M u{\rm d}v_g=0$ and $\|\nabla_g u\|_{L^n(M)}\leq 1$ from
\cite{Fontana} by $\|u\|_{0,\tau}\leq 1$ for every $\tau>0.$ On the
other hand, when $(M,g)$ is an $n-$dimensional {\it compact}
Riemannian manifold {\it with smooth boundary} $\partial M$, then
Cherrier \cite{Cherrier} proved that for every $0\leq
\alpha<2^\frac{1}{1-n}\alpha_n$, 
\begin{equation}\label{Cherrier-egyenlotlenseg}
    \sup_{\displaystyle\int_M u{\rm d}v_g=0,\ \|\nabla_g u\|_{L^n(M)}\leq
1}\int_M e^{\alpha |u|^\frac{n}{n-1}}{\rm d}v_g<\infty,
\end{equation}
and the above constant  is sharp, i.e., if
$\alpha>2^\frac{1}{1-n}\alpha_n$, then the supremum in
(\ref{Cherrier-egyenlotlenseg}) is infinite.

The study of Moser-Trudinger inequalities  on {\it non-compact}
Riemannian manifolds is more delicate, the {\it curvature}
playing a crucial role. 
On one hand, Yang \cite[Theorem 2.3]{Yang-JFA2012} proved that if
$(M,g)$ is an $n-$dimensional complete non-compact Riemannian
manifold with Ricci curvature bounded from below and positive
injectivity radius, then for every $\alpha\in [0,\alpha_n)$, there
exists   $\tau>0$ such that $({\bf MT})_{\alpha,\tau}^1$ holds on
$(M,g),$ while for every $\alpha>\alpha_n$ and $\tau>0$, $({\bf
MT})_{\alpha,\tau}^1$ fails on $(M,g).$ We emphasize that Yang's
result deeply exploits the existence of lower bounds on the harmonic
radius in terms of bounds on the Ricci curvature and the injectivity
radius, see Hebey \cite[Theorems 1.2 \& 1.3]{Hebey}.  On the other
hand, by using the arguments from Lam and Lu \cite{LamLu-Adv} and
fine estimates on the density function of the volume form, Yang, Su
and Kong \cite{YSK} proved that $({\bf MT})_{\alpha,\tau}^0$ holds
on every Hadamard manifold $(M,g)$ for every $\alpha\in
[0,\alpha_n]$ and $\tau>0,$ and $\alpha_n$ is again sharp;
furthermore, as a consequence of Yang \cite[Proposition
2.1]{Yang-JFA2012},  the embedding $W^{1,n}(M)\hookrightarrow
L^p(M)$ is continuous for every $p\in [n,\infty).$

In the next chapter we shall state and comment our achievements;
first, some theoretical results are established and then we present
an application on homogeneous Hadamard manifolds.


\section{Statement of main results}\label{sect1.2}

\subsection{Theoretical results:  validity of Moser-Trudinger inequalities} A first statement concerns the failure of
Moser-Trudinger inequalities in two different settings without any
curvature restriction.

\begin{proposition}\label{prop-veges-terfogat} 
Let $(M,g)$ be an $n-$dimensional  complete  Riemannian manifold,
$n\geq 2$. The following statements hold$:$
\begin{itemize}
    \item [{\rm (i)}] If $(M,g)$ is non-compact with ${\rm
    Vol}_g(M)<\infty$ then for any $\alpha>0$ and $\tau>0$, the Moser-Trudinger
inequalities $({\bf MT})_{\alpha,\tau}^i$ fail on $(M,g)$, $i\in
\{0,1\};$
 \item[{\rm (ii)}] For any $\alpha>\alpha_n$ and $\tau>0$,  the inequalities $({\bf MT})_{\alpha,\tau}^i$ fail on $(M,g)$, $i\in
\{0,1\}.$
\end{itemize}
\end{proposition}

According to Proposition \ref{prop-veges-terfogat}, Moser-Trudinger
inequalities $({\bf MT})_{\alpha,\tau}^i$ on any $n-$dimensional
non-compact complete Riemannian manifold $(M,g)$ are relevant
whenever ${\rm
    Vol}_g(M)=\infty$ and the parameter $\alpha$ belongs
to the subcritical interval
  $[0,\alpha_n].$

    Let $(M,g)$ be an
$n-$dimensional complete Riemannian manifold and $\Omega$ be a
smooth open subset in $M$, $n\geq 2$. We define the {\it
$n-$isoperimetric constant of $\Omega$} as
$$\mathcal I_n(\Omega,g):=\inf_{A}\frac{{\rm Area}_g(\partial A)}{{\rm Vol}_g(A)^{1-\frac{1}{n}}},$$
where $A$ varies over open sets of $\Omega$ having compact closure
and smooth boundary.  Hereafter, Area$_g(\partial A)$ stands for the
area of $\partial A$ with respect to the metric induced on $\partial
A$ by $g$, and Vol$_g(A)$ is the volume of $A$ with respect to $g$.
By considering geodesic balls $A:=B_x(r)$ in $\Omega\subset M$ with
$r\to 0^+,$ one clearly has
\begin{equation}\label{isop-eeee}
    \mathcal I_n(\Omega,g)\leq n\omega_n^\frac{1}{n},
\end{equation}
the number $n\omega_n^\frac{1}{n}$ being the $n-$dimensional
Euclidean isoperimetric ratio. For later use, let
\begin{equation}\label{isoperimetric-constant}
   {\rm Isop}{(\Omega,g)}:=\frac{\mathcal
   I_n(\Omega,g)}{n\omega_n^\frac{1}{n}}\in [0,1]
\end{equation}
be the {\it normalized $n-$isoperimetric constant of $\Omega$}.

 By using rearrangement
arguments on Riemannian manifolds in the spirit of Aubin-Hebey, see
\cite{Aubin, Hebey}, we prove the following quantitative result
which states a connection between the isoperimetric data of an open
set $\Omega\subset M$ and Moser-Trudinger inequalities on
$(\Omega,g)$:

\begin{lemma}\label{lemma-isop} 
Let $(M,g)$ be an $n-$dimensional  complete  Riemannian manifold,
$n\geq 2$, and $\Omega$ be a smooth open subset in $M$ such that
${\rm Isop}{(\Omega,g)}>0$. The following statements hold:
\begin{itemize}
  \item[{\rm (i)}] If\ ${\rm Vol}_g(\Omega)<\infty$, for  $\alpha\in \left[0,{\rm Isop}{(\Omega,g)}^\frac{n}{n-1}\alpha_n\right]$ and $u\in
C_0^\infty(\Omega)$ with $\|\nabla_g u\|_{L^n(\Omega)}\leq 1,$ one has
$$\int_\Omega \Phi_n(\alpha|u|^\frac{n}{n-1}){\rm d}v_g\leq M_0\|\nabla_g u\|_{L^n(\Omega)}^n{\rm Vol}_g(\Omega),$$
where $M_0>0$ is from {\rm (\ref{Moser-sharp}).}
    \item [{\rm (ii)}] For any $\tau>0$ and $\alpha\in \left[0,\min\left\{\tau^\frac{n}{n-1}, {\rm
    Isop}{(\Omega,g)}^\frac{n}{n-1}\right\}\alpha_n\right]$, the Moser-Trudinger
inequalities $({\bf MT})_{\alpha,\tau}^i$ hold on $(\Omega,g)$,
$i\in \{0,1\}.$
\end{itemize}

\end{lemma}

It is worth to point out the consistency of Proposition
\ref{prop-veges-terfogat} (i) and Lemma \ref{lemma-isop},
respectively. Indeed, when $(M, g)$ is of finite volume then ${\rm
Isop}(M, g) = 0$; the latter fact can be checked  by taking the
test-sets $A:=M\setminus \overline B_x(r)$ and letting $r\to 0.$

 By exploring Lemma \ref{lemma-isop} (i) and Gromov's
covering lemma, we may characterize the validity of Moser-Trudinger
inequalities on  manifolds with Ricci curvature bounded from below in the spirit of Coulhon, Saloff-Coste \cite{CSC} and Varopoulos \cite{Var}:

\begin{theorem}\label{thm-Ricci} 
Let  $(M,g)$ be an $n-$di\-men\-sional  complete non-compact
Riemannian manifold $(n\geq 2)$ with Ricci curvature bounded from
below. 
Then the following statements are equivalent:
\begin{itemize}
  \item[{\rm (i)}] There exists $\alpha\in (0,\alpha_n]$ and
  $\tau>0$ such that $({\bf MT})_{\alpha,\tau}^0$ holds on $(M,g);$
  \item[{\rm (ii)}] There exists $\alpha\in (0,\alpha_n]$ and
  $\tau>0$ such that $({\bf MT})_{\alpha,\tau}^1$ holds on $(M,g);$
  \item[{\rm (iii)}]  $\inf_{x\in M}{\rm Vol}_g(B_x(1))>0.$
\end{itemize}
Moreover, any of the above statements imply that the embedding
$W^{1,n}(M)\hookrightarrow L^p(M)$ is continuous for every $p\in
[n,\infty).$
\end{theorem}

\begin{remark}\rm
(a) If $(M,g)$ is an $n-$dimensional complete non-compact Riemannian
manifold with Ricci curvature bounded from below and positive
injectivity radius,  it follows by Croke \cite{Croke-injectivity}
that $\inf_{x\in M}{\rm Vol}_g(B_x(1))>0.$ Therefore, we may apply
Theorem \ref{thm-Ricci} in order to prove the validity of $({\bf
MT})_{\alpha,\tau}^1$ on $(M,g)$ for {\it some} $\alpha\in (0,
\alpha_n]$ and
 $\tau>0,$ recovering partially the result of
Yang \cite[Theorem 2.3]{Yang-JFA2012}. Note that in Yang's result
the {\it positivity of the injectivity radius} is indispensable.
Furthermore, our argument shows that once the normalized
$n-$isoperimetric constant ${\rm Isop}{(M,g)}$ is close to 1,  the
value $\alpha$ for which $({\bf MT})_{\alpha,\tau}^1$ holds on
$(M,g)$ approaches the critical exponent $\alpha_n$, see Remark
\ref{closer-and-closer}.

(b) Following the approach from Carron \cite{Carron} and Hebey
\cite[Lemma 2.2]{Hebey}, we stress that the implication
(ii)$\Rightarrow$(iii) is valid on {\it generic} Riemannian
manifolds, see Yang \cite[Proposition 2.1]{Yang-JFA2012}.   A
similar argument
also works  for (i)$\Rightarrow$(iii). 
\end{remark}

A remarkable consequence of Theorem \ref{thm-Ricci} is as follows:

\begin{corollary}\label{2-dim-corollary} 
Let $(M,g)$ be a two-di\-men\-sional complete non-compact Riemannian
manifold with non-negative sectional curvature. Then there exists
$\alpha\in (0,4\pi]$ and
  $\tau>0$ such that $({\bf MT})_{\alpha,\tau}^i$ holds on $(M,g),$
$i\in \{0,1\}$.
\end{corollary}

\begin{remark}\rm
(a) On one hand, Corollary \ref{2-dim-corollary} {\it cannot} be
deduced from Yang \cite[Theorem 2.3]{Yang-JFA2012} since no lower
bound for the injectivity radius can be guaranteed. Indeed, Croke
and Karcher \cite{CK} modified the paraboloid of revolution by
gluing to it a sequence of disjoint tangential cones in order to
obtain a hypersurface with positive sectional curvature and {\it
zero injectivity radius}. On the other hand, under the assumptions
of Corollary \ref{2-dim-corollary} it follows by \cite[Theorem
A]{CK} that $${\rm Vol}_g(B_x(r))\geq C_Mr^2\ \ {\rm for\ every}\
x\in M \ {\rm and}\ 0\leq r\leq 1,$$ the constant $C_M\in (0,\pi]$
depending only on $(M,g)$. Thus, it remains to apply Theorem
\ref{thm-Ricci} to conclude the proof of Corollary
\ref{2-dim-corollary}.

(b) We emphasize that Corollary \ref{2-dim-corollary} is {\it sharp}
with respect to the {\it dimension}. Indeed, one can construct
convex hypersurfaces $H$ in $\mathbb R^{n+1}$ with $n\geq 3$, having
positive sectional curvature and $\inf_{x\in H}{\rm
Vol}_g(B_x(1))=0$, see Croke and Karcher \cite[p. 755]{CK}.
Consequently, by Theorem \ref{thm-Ricci} and Proposition
\ref{prop-veges-terfogat} (ii),  the Moser-Trudinger inequalities
$({\bf MT})_{\alpha,\tau}^i$ fail on $H$ for every $\alpha>0$,
  $\tau>0$, and $i\in \{0,1\}$.
\end{remark}

Another byproduct of Lemma \ref{lemma-isop}  is a sharp
Moser-Trudinger inequality on Hadamard manifolds (simply connected, complete Riemannian manifold with non-positive sectional curvature):

\begin{corollary}\label{cor-Cartan Hadamard}
Let $(M,g)$ be an $n-$di\-men\-sional  Hadamard manifold  $(n\geq
2)$ which satisfies the Cartan-Hadamard conjecture. Then for every
$\alpha\in
  \left[0,\alpha_n\right]$ and $\tau\geq 1$, the Moser-Trudinger
inequalities $({\bf MT})_{\alpha,\tau}^i$ hold on $(M,g)$, $i\in
\{0,1\}.$ Moreover, the embedding $W^{1,n}(M)\hookrightarrow L^p(M)$
is continuous for every $p\in [n,\infty).$
\end{corollary}

\begin{remark}\rm \label{remark-conjecture}
Given an $n-$dimen\-sional Ha\-da\-mard manifold, $n\geq 2$, the
Cartan-Hadamard conjecture is equivalent to $ {\rm Isop}{(M,g)}=1$,
i.e., for every bounded open set $A\subset M$ with smooth boundary,
one has
\begin{equation}\label{cartan-hadamard-conjecture}
    {\rm Area}_g(\partial A) \geq n\omega_n^\frac{1}{n}  {\rm
    Vol}_g(A)^\frac{n-1}{n},
\end{equation}
see Aubin \cite{Aubin}. Note that the Cartan-Hadamard conjecture
holds on any Hadamard manifold of dimension
  $2$, cf. Beckenbach and Rad\'o \cite{BR} and  Weil \cite{Weil}, of dimension
$3$, cf. Kleiner \cite{Kleiner}, and of dimension $4$, cf. Croke
\cite{Croke}. We also notice that Corollary \ref{cor-Cartan Hadamard}
has been proved in \cite[Theorem 1.2]{YSK} without requiring the
validity of the Cartan-Hadamard conjecture;  the approach in
\cite{YSK} is based on fine estimates for Jacobi fields.
\end{remark}

\noindent


\subsection{Application: elliptic PDE on Hadamard manifolds with critical nonlinearity}\label{sect1.3} In the sequel, we shall present an application of the sharp Moser-Trudinger inequalities (Corollary \ref{cor-Cartan Hadamard} and \cite[Theorem 1.2]{YSK}) by considering the model
 elliptic problem
$$-\Delta_{n,g} u + |u|^{n-2}u=
f(u)\ {\rm in}\ M,\eqno{(\mathcal P)}$$
%
where $n\geq 2,$ $(M,g)$ is an $n-$dimensional Hadamard manifold,
$\Delta_{n,g}u={\rm div}_g(|\nabla_g u|^{n-2}\nabla_g u)$ is the
$n-$Laplace-Beltrami operator on $(M,g)$, and the continuous
function $f:[0,\infty)\to \mathbb R$ satisfies the following
hypotheses:
\begin{itemize}
  \item[$(f_0)$] there exists $\gamma>n$ with
$f(s)=\mathcal{O}(s^{\gamma-1})$ as $s\to 0^+;$
  \item[$(f_1)$] there is $\alpha_0>0$ with
$f(s)=\mathcal{O}(\Phi_n(\alpha_0s^\frac{n}{n-1}))$ as $s\to
\infty,$ and
$\lim_{s\to \infty}sf(s)e^{-\alpha_0s^\frac{n}{n-1}}=\infty;$
\vspace{-0.2cm}
  \item[$(f_2)$] there exists $\mu>n$ such that
$0<\mu F(s)\leq sf(s)$ for every $s>0,$ where $F(s)=\int_0^sf(t){\rm
d}t;$
  \item[$(f_3)$] there exist $R_0>0$ and $A_0>0$ such that $F(s)\leq A_0
  f(s)$ for every $s\geq R_0.$
\end{itemize}

\begin{remark}\rm Let $n=2$ and $f:[0,\infty)\to \mathbb R$ be defined by
$f(s)=\min\{1,s\}(e^{s^2}-1).$ Then  $f$ satisfies hypotheses
$(f_0)-(f_3)$ with $\gamma=\mu=3$ and $\alpha_0=R_0=A_0=1$.
\end{remark}

 Let ${\rm Isom}_g(M)$ be
the group of isometries of $(M,g)$ and $G$ be a subgroup of ${\rm
Isom}_g(M)$. The {\it orbit of} $x\in M$ under the action of $G$ is
$O_G^x=\{\sigma (x):\sigma\in G\}$. A function $u:M\to \mathbb R$ is
$G-${\it invariant} if $u(\sigma(x))=u(x)$ for every $x\in M$ and
$\sigma\in G,$ i.e., $u$ is constant on the orbit $O_G^x.$
The {\it fixed point set of}  $G$ on $M$ is given by ${\rm
Fix}_M(G)=\{x\in M:\sigma (x)=x\ {\rm for\ all}\ \sigma\in G\}.$

We shall prove the following result:

\begin{theorem}\label{theorem-existence}
Let $(M,g)$ be an $n-$dimensional homogeneous Hadamard manifold
$(n\geq 2)$, and let $G$ be a compact connected subgroup of ${\rm
Isom}_g(M)$ such that ${\rm Fix}_M(G)=\{x_0\}$ for some $x_0\in M$
and  ${\rm Card}(O^x_G)=\infty$ for every $x\in M\setminus \{x_0\}$.
If $f:[0,\infty)\to \mathbb R$ satisfies hypotheses $(f_0)-(f_3)$,
then problem $(\mathcal P)$ has a non-zero, non-negative,
$G-$invariant weak solution.
\end{theorem}


\begin{remark}\rm

(i) A similar result to Theorem \ref{theorem-existence} has been established on $\mathbb R^2$ by de Figueiredo,   Miyagaki and Ruf \cite{FMR}.   The novelty of Theorem \ref{theorem-existence} is twofold.
First,  no restriction is imposed on the boundedness from below of
the Ricci curvature on $(M,g)$ as in Yang \cite[Theorem
2.7]{Yang-JFA2012}. Second, Theorem \ref{theorem-existence} seems to
be the first existence result on non-compact Riemannian manifolds
involving exponential terms, by exploring deep features of the
isometric group in order to overcome some compactness. In order to
recover the non-compactness of the space (even in the Euclidean
case), instead of the left-hand side of $(\mathcal P)$, most of the
authors considered operators of the form $u\mapsto -\Delta_{n,g} u +
V(x)|u|^{n-2}u$ where $V$ is {\it coercive,} i.e., $V(x)\to \infty$
as $d_g(x_0,x)\to \infty$ for some $x_0\in M$ fixed, see e.g.
Adimurthi and  Yang \cite{Adimurthi-Yang}, do \'O \cite{doO}, do \'O
and Yang \cite{doO-Yang}, Lam and Lu \cite{LamLu}, Yang
\cite{Yang-JFA2012}. Under this coercivity assumption a
Rabinowitz-type argument shows that the weighted Sobolev space
$W^{1,n}_V(M)=\{u\in W^{1,n}(M):\displaystyle\int_M V(x)|u|^n{\rm d}v_g<\infty\}$
 is compactly embedded into $L^p(M)$, $p\in [n,\infty)$. In our
 case such approach fails. However, in order to prove Theorem \ref{theorem-existence}, we shall combine the
 principle of symmetric criticality of Palais \cite{Palais} with a
 recent characterization of compactness of invariant Sobolev spaces \`a la Lions (see \cite{Lions}) under the action of isometries, see Skrzypczak and Tintarev
 \cite{S-Tintarev}. As far as we know, the only result
 for $V\equiv 1$ in $\mathbb R^n$ has been provided  recently
 by do \'O, de Souza, de Medeiros and Severo \cite{doO-JDE-2014} via a Lions-type concentration-compactness argument.

(ii) Theorem \ref{theorem-existence} is new even in the Euclidean
case where one can choose certain subgroups $G$ of the special
orthogonal group in $\mathbb R^n.$
  Further examples will be provided in \S \ref{sect-4} on
the $n-$dimensional hyperbolic space, and on the open convex cone of
symmetric positive definite matrices endowed with a trace-type
scalar product.

\end{remark}


\section{Preliminaries}
\subsection{Generic notions} In order the paper to be self-contained, we recall those ingredients
from Riemannian geometry which will be used throughout the paper.
Let $(M, g)$ be an $n-$dimensional Riemannian manifold, $T_xM$ be
the tangent space at $x\in M,$ $TM=\cup_{x\in M}T_x M$ be the
tangent bundle, and $d_g:M\times M\to [0,\infty)$ be the induced
metric function by the Riemannian metric $g$. As usual,
  let $B_x(r)=\{y\in M:d_g(x,y)<r\}$ and
$\overline B_x(r)=\{y\in M:d_g(x,y)\leq r\}$ be the open and closed
geodesic balls with center $x\in M$ and radius $r>0$, respectively.
If
 ${\text d}v_g$ is the canonical volume
element on $(M,g)$, the volume of an open bounded set $\Omega\subset
M$ is Vol$_g(\Omega)=\displaystyle\int_\Omega {\text d}v_g=\mathcal H^n(\Omega)$,
where $\mathcal H^n(S)$ denotes the $n-$dimensional Hausdorff
measure of $\Omega$ with respect to the metric $d_g$. Let ${\text
d}\sigma_g$ be the $(n-1)-$dimensional Riemann measure induced on
$\partial \Omega$ by $g$; then
 Area$_g(\partial \Omega)=\displaystyle\int_{\partial \Omega} {\text d}\sigma_g=\mathcal
H^{n-1}(\partial \Omega)$ is the area of $\partial \Omega$ with
respect to the metric $g$.  For further use, $B_0(\delta)$, ${\rm
d}x$, ${\rm d}\sigma_e$, Vol$_e(S)$ and Area$_e(S)$ denote the
Euclidean counterparts of the above notions when  $S\subset \mathbb
R^n$.

The behavior of the volume of small geodesic balls can be expressed
as follows; for every $x\in M$ we have (see Gallot, Hulin and
Lafontaine \cite[Theorem 3.98]{GHL}):
\begin{equation}\label{volume-aszimpt}
    {\rm Vol}_g(B_x(\rho))=\omega_n \rho^n(1+o(\rho))\ {\rm as}\ \rho\to
    0.
\end{equation}

The manifold $(M,g)$ has {\it Ricci curvature bounded from below} if
there exists $k\in \mathbb R$ such that ${\rm Rc}_{(M,g)} \geq kg$
in the sense of bilinear forms, i.e., ${\rm Rc}_{(M,g)}(X,X) \geq
k|X|_x^2$ for every $X\in T_xM$ and $x\in M,$ where ${\rm
Rc}_{(M,g)}$ is the Ricci curvature, and $|X|_x$ denotes the norm of
$X$ with respect to the metric $g$ at the point $x$. For simplicity
of notation, $\langle\cdot,\cdot\rangle_x$ denotes the scalar
product $g_x$ on $T_x M$ induced by the metric $g$. When no
confusion arises, if $X,Y\in T_x M$, we simply write $|X|$ and
$\langle X,Y\rangle$ instead of $|X|_x$ and $\langle X,Y\rangle_x$,
respectively.

In the sequel, $V_k(\rho)$ shall denote the volume of a ball of
radius $\rho$ in the $n-$dimensional simply connected, complete
Riemannian manifold of constant sectional curvature $k\in \mathbb
R.$ The behavior of the volume of large geodesic balls is given by
Bishop-Gromov and Bishop-Gunther:

\begin{proposition}\label{comparison} {\rm \cite[Theorem 3.101]{GHL}} Let $(M, g)$ be an  $n-$dimensional  complete
Riemannian manifold. The following statements hold:
\begin{itemize}
  \item[{\rm (i)}] If ${\rm Rc}_{(M,g)} \geq k(n-1)g$ for some $k\in
\mathbb R$, then $\rho\mapsto \frac{{\rm
Vol}_g(B_x(\rho))}{V_k(\rho)}$ is non-increasing for every $x\in M$.
In particular, by {\rm (\ref{volume-aszimpt})}, one has ${\rm
Vol}_g(B_x(\rho))\leq V_k(\rho)$ for every $\rho\geq 0$ and $x\in
M.$
  \item[{\rm (ii)}] If the sectional curvature of $(M,g)$ is bounded
  from above by $k\in \mathbb R$, then ${\rm
Vol}_g(B_x(\rho))\geq V_k(\rho)$ for every $\rho\geq 0$ and $x\in
M.$
\end{itemize}
\end{proposition}

 The following result, which is a {\it local}
isoperimetric inequality on Riemannian manifolds with Ricci
curvature bounded from below, plays a crucial role in the proof of
Theorem \ref{thm-Ricci}.

\begin{proposition}\label{prop-izometrikus} {\rm \cite[Lemma 3.2]{Hebey}} Let $(M, g)$ be an  $n-$dimensional  complete Riemannian
manifold whose Ricci curvature satisfies ${\rm Rc}_{(M,g)} \geq kg$
for some $k\in \mathbb R$, and suppose that there exists $v > 0$
such that ${\rm Vol}_g(B_x(1))\geq v$ for every $x\in M$. Then there
exist two positive constants $C_0 = C(n, k, v)$ and $\eta_0 =
\eta(n, k, v),$ depending only on $n, k,$ and $v,$ such that for any
open set $\Omega\subset M$ with smooth boundary and compact closure,
if ${\rm Vol}_g(\Omega)\leq \eta_0$, then $${\rm Area}_g(\partial
\Omega)\geq C_0 {\rm Vol}_g(\Omega)^\frac{n-1}{n}.$$
\end{proposition}

Gromov's covering lemma, whose proof is based on Proposition
\ref{comparison} (i), reads as follows:

\begin{proposition}\label{Gromov-partition} {\rm \cite[Lemma
1.1]{Hebey}} Let $(M, g)$ be an  $n-$dimensional  complete
Riemannian manifold whose Ricci curvature satisfies ${\rm
Rc}_{(M,g)} \geq kg$ for some $k\in \mathbb R$, and let $\rho>0$ be
fixed. Then there exists a sequence  $\{x_j\}_{j\in I}\subset M$
$($with $I$ countable$)$ such that for every $r\geq \rho:$
\begin{itemize}
  \item[{\rm (i)}] the family of sets $\{B_{x_j}(r)\}$ is a uniformly locally
  finite covering of $M$ and there exists an upper bound $N_0$ for
  this covering
  in terms of $n$, $\rho,$ $r$ and $k;$
  \item[{\rm (ii)}]  $B_{x_i}(\frac{\rho}{2})\cap B_{x_j}(\frac{\rho}{2})=\emptyset$ for
  every $i\neq j.$
\end{itemize}
\end{proposition}

\subsection{Compactness vs. symmetrization on Hadamard manifolds}
Let $p\in [ 1,\infty).$ The norm of $L^p(M)$ is given by
$\|u\|_{L^p(M)}=\left(\displaystyle\int_M |u|^p{\rm
d}v_g\right)^\frac{1}{p}$ while $\|\cdot\|_{L^\infty(M)}$ denotes
the usual supremum-norm. Let $u:M\to \mathbb R$ be a function of
class $C^1.$ If $(x^i)$ denotes the local coordinate system on a
coordinate neighborhood of $x\in M$, and the local components of the
differential of $u$ are  $u_i=\frac{\partial u}{\partial x^i}$, then
the local components of the gradient  $\nabla_g u$ are
$u^i=g^{ij}u_j$. Here, $g^{ij}$ are the local components of
$g^{-1}=(g_{ij})^{-1}$. In particular, for every $x_0\in M$ one has
\begin{equation}\label{dist-gradient}
    |\nabla_g d_g(x_0,\cdot)|=1\ {\rm  a.e.\ on}\
M.
\end{equation}
The $L^n(M)$ norm of  $\nabla_g u(x)\in T_xM$ is given by
$$\|\nabla_g u\|_{L^n(M)}=\left(\displaystyle\int_M |\nabla_gu|^n{\rm
d}v_g\right)^\frac{1}{n},$$ while the space $W^{1,n}(M)$ is the
completion of $C_0^\infty(M)$ with respect to the norm
$\|\cdot\|_{0,1}.$

 In the sequel we adapt the main results from Skrzypczak
and Tintarev \cite{S-Tintarev} to our setting concerning the Sobolev
spaces in the presence of group-symmetries; for a similar approach
see also Hebey and Vaugon \cite{HV}. When $(M,g)$ is a Hadamard
manifold, the embedding $W^{1,n}(M)\hookrightarrow L^p(M)$ is
continuous for every $p\in [n,\infty)$ (cf. Corollary \ref{cor-Cartan
Hadamard}), but not compact. By exploiting the fact that the
embedding $W^{1,n}(M)\hookrightarrow L^p(M)$ is (weakly) cocompact
relative to the isometry group ${\rm Isom}_g(M)$ for every $p\in
(n,\infty)$, one can state the following result:

\begin{proposition}\label{Tintarev} {\rm \cite[Theorem 1.3 \& Proposition 3.1]{S-Tintarev}} Let $(M, g)$ be an  $n-$dimensional homogeneous Ha\-da\-mard
manifold and $G$ be a compact connected subgroup of ${\rm
Isom}_g(M)$ such that ${\rm Fix}_M(G)$ is a singleton. Then the
subspace of $G-$invariant functions of $W^{1,n}(M)$, i.e.,
$$W_G^{1,n}(M)=\{u\in W^{1,n}(M): u\circ \sigma=u\ for\ all \
\sigma\in G\}$$ is compactly embedded into $L^p(M)$ for every $p\in
(n,\infty)$.
\end{proposition}



We conclude this section with the principle of symmetric criticality
of Palais \cite{Palais}. A group $G$ {\it acts continuously} on a
real Banach space $W$ by an application $[\sigma,u]\mapsto \sigma u$
from $G\times W$ to $W$ if this map itself is continuous on $G\times
W$ and
\begin{itemize}
  \item $\sigma_{\rm id}u=u$ for every $u\in W$, where $\sigma_{\rm id}\in
  G$ is the identity element of $G$;
  \item $(\sigma_1\sigma_2)u=\sigma_1(\sigma_2u)$ for
every $\sigma_1,\sigma_2\in G$ and $u\in W$;
  \item $u\mapsto \sigma u$
is linear for every $\sigma\in G$.
\end{itemize}

\begin{proposition}\label{Palais-PSC} {\rm \cite{Palais}}
Let $W$ be a real Banach space, $G$ be a compact topological group
acting continuously on $W$ by a map $[\sigma,u]\mapsto \sigma u$
from $G\times W$ to $W$, and $h:W\to \mathbb R$ be a $G-$invariant
$C^1-$function, i.e., $h(\sigma u)=h(u)$ for every $(\sigma,u)\in
G\times W$. If $u_G\in {\rm Fix}_W(G)=\{u\in W:\sigma u=u\ for\ all
\ \sigma\in G\}$ is a critical point of $h_G=h|_{{\rm Fix}_W(G)},$
then $u_G$ is also a critical point of $h.$
\end{proposition}

\section{Proof of theoretical results}

{\it Proof of Proposition \ref{prop-veges-terfogat}.} (i) It is
enough to prove the statement for $({\bf MT})_{\alpha,\tau}^1$. By
contradiction, let us assume that 
$({\bf MT})_{\alpha,\tau}^1$ holds on $(M,g)$ for some $\alpha>0$
and $\tau>0$. Due to Yang \cite[Proposition 2.1]{Yang-JFA2012},
there is $v>0$ such that
 ${\rm Vol}_g(B_x(1))\geq v$ for every $x\in M.$ A similar
 argument as in Hebey \cite[pp. 53-54]{Hebey} based on Zorn lemma shows that there
 exists a sequence
 $\{x_i\}_{i\in I}\subset M$ such that $B_{x_i}(1)\cap B_{x_j}(1)=\emptyset$  for every $i\neq j$ and $M=\bigcup_{i\in I}B_{x_i}(2).$
 Note that $$+\infty> {\rm Vol}_g(M)\geq \sum_{i\in I}{\rm Vol}_g(B_{x_i}(1))\geq {\rm
 Card}(I)v.$$ Therefore, $I$ is finite, which implies together
 with the Hopf-Rinow theorem that $M$ is covered by a finite number of relatively compact sets; thus  $M$ is compact, a contradiction.

(ii) A similar statement for $({\bf MT})_{\alpha,\tau}^1$ is presented
by Yang \cite{Yang-JFA2012} on Riemannian manifolds with Ricci
curvature bounded from below and positive injectivity radius. In
fact, since the Moser-type truncation functions are locally
constructed, the proof in \cite{Yang-JFA2012} works in generic
Riemannian manifolds as well; for completeness we provide its
 proof since some parts will be used later on.

Let $x_0\in M$ be arbitrarily fixed and denote by $i_{x_0}$ the
injectivity radius at $x_0;$ clearly, $i_{x_0}>0.$ Choose also
$\varepsilon_0\in (0,i_{x_0})$ sufficiently small such that it
belongs to the range of (\ref{volume-aszimpt}). For every
$\varepsilon\in (0,\varepsilon_0)$, we introduce the Moser-type
truncation function $u_\varepsilon:M\to [0,\infty)$ defined by
\begin{equation}\label{Moser-function-1}
    u_\varepsilon(x)=\min\left\{\left(\frac{\log\frac{\varepsilon_0}{d_g(x_0,x)}}{\log\frac{\varepsilon_0}{\varepsilon}}\right)_+,1\right\},
\end{equation}
where $r_+=\max\{0,r\}$ for $r\in \mathbb R.$ The functions
$u_\varepsilon$ can be approximated by elements from $C_0^\infty(M)$
and we shall see that $u_\varepsilon \in W^{1,n}(M).$
Indeed, on one hand, a simple computation combined with
(\ref{volume-aszimpt}), (\ref{dist-gradient})
 and the layer cake representation gives
\begin{eqnarray*}
  \|\nabla_g u_\varepsilon\|_{L^n(M)}^n 
   &=&n\omega_n{\left(\log\frac{\varepsilon_0}{\varepsilon}\right)^{1-n}}\left(1+\mathcal{O}\left(\left(\log\frac{\varepsilon_0}{\varepsilon}\right)^{-1}\right)\right),
\end{eqnarray*}
as $\varepsilon\to 0.$ On the other hand, again by the layer cake
representation and (\ref{volume-aszimpt}), we obtain that
\begin{eqnarray*}
 \| u_\varepsilon\|_{L^n(M)}^n 
   &=&\mathcal{O}\left(\left(\log\frac{\varepsilon_0}{\varepsilon}\right)^{-n}\right),
\end{eqnarray*}
as $\varepsilon\to 0.$ Consequently, since $n\omega_n=\omega_{n-1}$,
if $\tau>0$ is arbitrarily fixed, one has that
\begin{equation}\label{ueps-becsles}
    \|u_\varepsilon\|_{1,\tau}=\omega_{n-1}^\frac{1}{n}{\left(\log\frac{\varepsilon_0}{\varepsilon}\right)^\frac{1-n}{n}}\left(1+\mathcal{O}\left(\left(\log\frac{\varepsilon_0}{\varepsilon}\right)^{-\frac{1}{n}}\right)\right).
\end{equation}
Therefore, if $\alpha>\alpha_n=n\omega_{n-1}^\frac{1}{n-1}$, by
relations (\ref{volume-aszimpt}) and (\ref{ueps-becsles}) it follows
that
\begin{eqnarray*}
  S_{\alpha,\tau}^1(M,g) &\geq& \lim_{\varepsilon\to 0}\int_M \Phi_n\left(\alpha\frac{u_\varepsilon^\frac{n}{n-1}}{\|u_\varepsilon\|_{1,\tau}^\frac{n}{n-1}}\right){\rm d}v_g
   \geq \lim_{\varepsilon\to 0}\int_{B_{x_0}(\varepsilon)} \Phi_n\left(\frac{\alpha}{\|u_\varepsilon\|_{1,\tau}^\frac{n}{n-1}}\right){\rm d}v_g\\
   &=&\lim_{\varepsilon\to 0}\left({\rm Vol}_g(B_{x_0}(\varepsilon))\Phi_n\left({\alpha}{\|u_\varepsilon\|_{1,\tau}^\frac{n}{1-n}}\right)\right)
   =\varepsilon_0^{\alpha \omega_{n-1}^\frac{1}{1-n}}\lim_{\varepsilon\to
   0}\varepsilon^{n-\alpha \omega_{n-1}^\frac{1}{1-n}}\\
   &=&+\infty,
\end{eqnarray*}
which means that the Moser-Trudinger inequality $({\bf
MT})_{\alpha,\tau}^1$ fails on $(M,g)$. Since
$S_{\alpha,\tau}^1(M,g)\leq S_{\alpha,\tau}^0(M,g)$, $({\bf
MT})_{\alpha,\tau}^0$ also fails on
$(M,g)$, which concludes the proof. \hfill $\square$\\

{\it Proof of Lemma \ref{lemma-isop}.} We divide the proof into four
steps.

{\it \underline{Step 1}: choice of test functions.} Since for every
$u\in W^{1,n}(M)$ we have $|\nabla_g u|=|\nabla_g |u||$ a.e. on $M,$
classical Morse theory and density argument show that
Moser-Trudinger inequalities on $(M,g)$ are sufficient to be
considered for continuous
 test functions $u:M\to [0,\infty)$ having compact support $S\subset
 M$, where $S$ is enough smooth, $u$ being of class
 $C^\infty$ in $S$ and having only non-degenerate critical points in $S.$

{\it  \underline{Step 2}: P\'olya-Szeg\H{o}-type inequality.}   Let $\Omega\subset M$ be an open set  
and we consider a non-negative function $u\in C_0^\infty(\Omega)$
with the properties from Step 1. To this function $u,$ we associate
its
 Euclidean rearrangement function $u^*:\mathbb R^n\to
 [0,\infty)$ which is
 radially symmetric, non-increasing in $|x|$, and for every $t>0$ is defined by
 \begin{equation}\label{vol-egyenloseg}
    {\rm Vol}_e(\{x\in \mathbb R^n:u^*(x)>t\})={\rm Vol}_g(\{x\in
    \Omega:u(x)>t\}).
 \end{equation}
It is clear by (\ref{vol-egyenloseg}) that
\begin{equation}\label{ket-terfogat}
  {\rm Vol}_g(\Omega)\geq  {\rm Vol}_g(S)={\rm Vol}_g({\rm supp}(u))={\rm Vol}_e({\rm
supp}(u^*))={\rm Vol}_e(B_0(R_S))
\end{equation}
 for some $R_S>0$. By the layer cake representation, for every $q>0$ one
has
\begin{eqnarray}\label{L^q-normak}
 \nonumber \int_\Omega u^q{\rm d}v_g &=& \int_0^\infty{\rm Vol}_g(\{x\in
    \Omega:u^q(x)>t\}){\rm d}t= \int_0^\infty{\rm Vol}_e(\{x\in \mathbb R^n:(u^*)^q(x)>t\}) \\
   &=& \int_{B_0(R_S)}(u^*)^q{\rm d}x.
\end{eqnarray}
For abbreviation, we consider the sets
$$A_t=\{x\in
    \Omega:u(x)>t\},\ \ \ A^*_t=\{x\in \mathbb R^n:u^*(x)>t\}$$
    for every $0<t<\|u\|_{L^\infty(M)}=\|u^*\|_{L^\infty(\mathbb R^n)}.$
The boundaries of $A_t$ and $A_t^*$ are exactly the
 level sets $\partial A_t=u^{-1}(t)\subset S\subset \Omega$ and $\partial A_t^*=(u^*)^{-1}(t)\subset \mathbb R^n,$ which are regular. Since $u^*$ is radially symmetric,
the set $\partial A_t^*$ is an $(n-1)-$dimensional sphere for every
$0<t<\|u\|_{L^\infty(M)}=\|u^*\|_{L^\infty(\mathbb R^n)}.$
Therefore, by relation (\ref{vol-egyenloseg}) we have that
\begin{eqnarray}\label{nocsak-ejsze}
 \nonumber  {\rm Area}_g(\partial A_t) &\geq &\mathcal
I_n(\Omega,g) {\rm Vol}_g(A_t)^\frac{n-1}{n} = \mathcal
I_n(\Omega,g) {\rm Vol}_e(A^*_t)^\frac{n-1}{n} \\
   &=& {\rm Isop}{(\Omega,g)}{\rm Area}_e(\partial A_t^*).
\end{eqnarray}
Let $$V(t)={\rm Vol}_g(A_t)={\rm Vol}_e(A^*_t).$$  The co-area
formula (see Chavel \cite[pp. 302-303]{Chavel}) gives
\begin{equation}\label{V-deriv}
    V'(t)=-\int_{\partial A_t}\frac{1}{|\nabla_g u|}{{\rm d}
\sigma}_g=-\int_{\partial A_t^*}\frac{1}{|\nabla u^*|}{{\rm d}
\sigma}_e.
\end{equation}
 Since $|\nabla u^*|$ is
constant on the sphere $\partial A_t^*$, by (\ref{V-deriv}) one has
that
\begin{equation}\label{V-DER-MASIK}
V'(t)=-\frac{{\rm Area}_e(\partial A_t^*)}{|\nabla u^*(x)|},\ x\in
\Gamma_t^*.
\end{equation}
By (\ref{V-deriv}) again and H\"older's inequality it turns out that
$${\rm Area}_g(\partial A_t)=\int_{\partial A_t}{{\rm d} \sigma}_g\leq\left(-V'(t)\right)^\frac{n-1}{n}\left(\int_{\partial A_t}{|\nabla_g u|^{n-1}}{{\rm
d} \sigma}_g\right)^\frac{1}{n}.$$ For every
$0<t<\|u\|_{L^\infty(M)}$, by using (\ref{nocsak-ejsze}) and
(\ref{V-DER-MASIK}), it follows that
\begin{eqnarray*}
  \int_{\partial A_t}{|\nabla_g u|^{n-1}}{{\rm
d} \sigma}_g &\geq& {\rm Area}_g(\partial A_t)^n \left(-V'(t)\right)^{1-n}\\
   &\geq & {\rm Isop}{(\Omega,g)}^n{\rm Area}_e(\partial A_t^*)^n \left(\frac{{\rm Area}_e(\partial A_t^*)}{|\nabla u^*(x)|}\right)^{1-n}\ \ \ \ \ \ \ \ (x\in \Gamma_t^*)\\
   &=& {\rm Isop}{(\Omega,g)}^n\int_{\partial
A_t^*}{|\nabla u^*|^{n-1}}{{\rm d} \sigma}_e.
\end{eqnarray*}
The co-area formula and the above estimate  give a
P\'olya-Szeg\H{o}-type inequality
\begin{eqnarray}\label{Polya-Szego}
  \nonumber \int_{\Omega}{|\nabla_g u|^{n}}{{\rm d}
}v_g &=&\nonumber \int_0^\infty\int_{\partial A_t}{|\nabla_g u|^{n-1}}{\rm d }\sigma_g {\rm d}t \\
   &\geq& \nonumber{\rm Isop}{(\Omega,g)}^n\int_0^\infty
\int_{\partial A_t^*}{|\nabla u^*|^{n-1}}{{\rm d} \sigma}_e ={\rm
Isop}{(\Omega,g)}^n\int_{\mathbb R^n}{|\nabla u^*|^{n}}{\rm
d}x\nonumber\\&=& {\rm Isop}{(\Omega,g)}^n\int_{B_0(R_S)}{|\nabla
u^*|^{n}}{\rm d}x.
\end{eqnarray}
Now, we are ready to prove the claims (i) and (ii).

{{\it  \underline{Step 3}}: {\it proof of} (i).}  Let $\Omega$ be an
open subset of $M$ such that ${\rm Isop}{(\Omega,g)}>0$,
Vol$_g(\Omega)<\infty$ and let $u\in C_0^\infty(\Omega)$ be a
non-negative and non-zero function with the properties from Step 1
and
\begin{equation}\label{kisebb-mint-egy}
    \|\nabla_g u\|_{L^n(\Omega)}\leq 1.
\end{equation}
Let $\tilde u=\frac{u}{\|\nabla_g u\|_{L^n(\Omega)}}.$ Applying the
arguments from Step 2 for the function $\tilde u$, by
(\ref{Polya-Szego}) it follows that
\begin{equation}\label{sub-unitar}
    1=\|\nabla_g \tilde u\|_{L^n(\Omega)}\geq {\rm Isop}{(\Omega,g)}\|\nabla \tilde u^*\|_{L^n(B_0(R_S))},
\end{equation}
where $\tilde u^*$ is the Euclidean rearrangement function of
$\tilde u.$ Thus, for every $\alpha\in \left[0,{\rm
Isop}{(\Omega,g)}^\frac{n}{n-1}\alpha_n\right]$, we have
\begin{eqnarray*}
  \int_\Omega \Phi_n\left(\alpha|u|^\frac{n}{n-1}\right){\rm d}v_g &=& \sum_{j=n-1}^\infty \frac{\alpha^j}{j!}\int_\Omega \tilde u^\frac{nj}{n-1}\|\nabla_g u\|_{L^n(\Omega)}^\frac{nj}{n-1}{\rm
  d}v_g\\&\leq& \|\nabla_g u\|_{L^n(\Omega)}^n\sum_{j=n-1}^\infty \frac{\alpha^j}{j!}\int_\Omega \tilde u^\frac{nj}{n-1}{\rm
  d}v_g\ \ \ \ \ \ \ \ \ \ \ \ \ \ \ \ \ \ \ \ \ \ \ \ \ {\rm [see\ (\ref{kisebb-mint-egy})]} \\ &=& \|\nabla_g u\|_{L^n(\Omega)}^n\sum_{j=n-1}^\infty \frac{\alpha^j}{j!}\int_{B_0(R_S)} (\tilde u^*)^\frac{nj}{n-1}{\rm
  d}x\ \ \ \ \ \ \ \ \ \ \ \ \ \ \ \ \ \  [{\rm see\ (\ref{L^q-normak}})]
\end{eqnarray*}

\begin{eqnarray*}
 \ \ \ \ \ \ \ \ \ \ \ \ \  \ \ \ \ \ \ \ \ \ \ \ \ \ \ \
&=&  \|\nabla_g u\|_{L^n(\Omega)}^n\int_{B_0(R_S)}
\Phi_n\left(\alpha(\tilde u^*)^\frac{n}{n-1}\right){\rm d}x\\
  &\leq&  \|\nabla_g u\|_{L^n(\Omega)}^n\int_{B_0(R_S)} \Phi_n\left(\alpha_n({\rm Isop}{(\Omega,g)}\tilde u^*)^\frac{n}{n-1}\right){\rm d}x\\&\leq& M_0 \|\nabla_g
u\|_{L^n(\Omega)}^n{\rm Vol}_e(B_0(R_S)) \ \ \ \ \ \ \ \ \ \ \ \ \ \
\ \ \ \ \  \ {\rm [see\ (\ref{sub-unitar})\ and\
(\ref{Moser-sharp})]}
\\&\leq& M_0 \|\nabla_g u\|_{L^n(\Omega)}^n{\rm Vol}_g(\Omega).\ \ \ \ \ \ \ \ \ \ \ \ \
\ \ \ \ \ \  \ \ \ \ \ \ \ \ \ \ \ \ \ \ \ \ \ \  {\rm [see\
(\ref{ket-terfogat})]}
\end{eqnarray*}

{{\it \underline{Step 4}}: {\it proof of} (ii).} Let us fix $\tau>0$
and $\alpha\in \left[0,\min\{\tau^\frac{n}{n-1}, {\rm
    Isop}{(\Omega,g)}^\frac{n}{n-1}\}\alpha_n\right]$, and  let $u\in
C_0^\infty(\Omega)$ be a non-negative and non-zero function with the
properties from Step 1 and $\|u\|_{0,\tau}\leq 1.$ Then, by
(\ref{Polya-Szego}) we have
\begin{eqnarray*}
  1 &\geq& \|u\|_{0,\tau}=\left(\int_{\Omega}(|\nabla_g u|^n+\tau^nu^n){\rm
d}v_g\right)^{1/n}\geq \left(\int_{\mathbb R^n}({\rm
Isop}{(\Omega,g)}^n|\nabla
u^*|^n+\tau^n(u^*)^n){\rm d}x\right)^{1/n} \\
   &\geq& \min\{{\rm
    Isop}{(\Omega,g)},\tau\}\|u^*\|_{0,1}.
\end{eqnarray*}
By this estimate, (\ref{L^q-normak}) and (\ref{Li-Ruf-eredmeny}) one
has
\begin{eqnarray*}
  \int_\Omega \Phi_n\left(\alpha|u|^\frac{n}{n-1}\right){\rm d}v_g &=& \int_{\mathbb R^n} \Phi_n\left(\alpha(u^*)^\frac{n}{n-1}\right){\rm d}x\\
   &\leq&  \int_{\mathbb R^n} \Phi_n\left(\alpha_n(\min\{{\rm
    Isop}{(\Omega,g)},\tau\} u^*)^\frac{n}{n-1}\right){\rm d}x\\&\leq&
    S_n^{LR},
\end{eqnarray*}
i.e., $({\bf MT})_{\alpha,\tau}^0$ holds on $(\Omega,g)$. This fact
also implies that $({\bf MT})_{\alpha,\tau}^1$ holds on
$(\Omega,g)$. \hfill $\square$\\

{\it Proof of Theorem \ref{thm-Ricci}.} (i)$\Leftrightarrow$(ii) is
trivial since the two norms $\|\cdot\|_{0,\tau}$ and
$\|\cdot\|_{1,\tau}$ are equivalent. More precisely, if $({\bf
MT})_{\alpha,\tau}^0$ holds on $(M,g)$ for some $\alpha>0$ and
$\tau>0,$ then $({\bf MT})_{\alpha,\tau}^1$ also holds on $(M,g)$;
conversely, if $({\bf MT})_{\alpha,\tau}^1$ holds on $(M,g)$ then
$({\bf MT})_{\tilde\alpha,\tau}^0$ holds on $(M,g)$ where $\tilde
\alpha=2^\frac{n}{1-n}\alpha$.

(ii)$\Rightarrow$(iii) is given in Yang \cite[Proposition
2.1]{Yang-JFA2012} for generic Riemannian manifolds; namely, if
$({\bf MT})_{\alpha,\tau}^1$ holds on $(M,g)$ for some $\alpha>0$
and $\tau>0$ then for every $q>n$, $x\in M$ and $r>0$, one has
$${\rm Vol}_g(B_x(r))\geq \min\left\{\frac{1}{2\tau Q},\frac{r}{2^\frac{2q-n}{q-n}Q}\right\}^\frac{nq}{q-n},$$
where $Q$ depends on $n,$ $q$, $\alpha$ and
$S_{\alpha,\tau}^1(M,g)<\infty.$

(iii)$\Rightarrow$(ii) Let $v > 0$ be such that ${\rm
Vol}_g(B_x(1))\geq v$ for every $x\in M$. According to Proposition
\ref{prop-izometrikus}, there exist two positive constants $C_0 =
C(n, k, v)$ and $\eta_0 = \eta(n, k, v),$ depending only on $n, k$
and $v,$ such that for any open set $\Omega\subset M$ with smooth
boundary and compact closure, if ${\rm Vol}_g(\Omega)\leq \eta_0$,
then  ${\rm Area}_g(\partial \Omega)\geq C_0 {\rm
Vol}_g(\Omega)^\frac{n-1}{n}.$ Consequently, one has
\begin{equation}\label{area-volume-333}
    {\rm Isop}{(\Omega,g)}\geq \frac{C_0}{n\omega_n^\frac{1}{n}}\ \ \ \ {\rm for\ all\ smooth\ open\ set}\ \Omega\subset M\ {\rm and}\ {\rm Vol}_g(\Omega)\leq \eta_0.
\end{equation}
Let $\rho_0>0$ be small enough such that
$\omega_n\rho_0^ne^{\sqrt{(n-1)|k|}\rho_0}\leq \eta_0.$
According to Proposition \ref{comparison} (i), it turns out that
\begin{equation}\label{volume-unif-estimate}
    {\rm Vol}_g(B_x(\rho_0))\leq \eta_0\ {\rm for\ all}\ x\in M.
\end{equation}
 By
Proposition \ref{Gromov-partition}, there exists a sequence
$\{x_j\}_{j\in \mathbb N}\subset M$ such that
$B_{x_i}(\frac{\rho_0}{4})\cap B_{x_j}(\frac{\rho_0}{4})=\emptyset$
for every $i\neq j$ and the family of geodesic balls
$B_{x_j}(\frac{\rho_0}{2})$ is a uniformly locally finite covering
of $M$, the number $N_0\in \mathbb N$ being the uniform upper bound
for this covering (which depends only on $\rho_0$, $n$ and $k$). Let
us fix $x_0\in M$ arbitrarily. For every $j\in \mathbb N,$ let
$$\psi_j(x)=\min\left\{\left(2-\frac{2}{\rho_0}d_g(x_j,x)\right)_+,1\right\},\ x\in M.$$
We have that $\psi_j\in W^{1,n}(M)$,
 $0\leq
\psi_j\leq 1$, $\psi_j(x)=1$ for every $x\in
B_{x_j}(\frac{\rho_0}{2}),$ $\psi_j(x)=0$ for every $x\in M\setminus
B_{x_j}(\rho_0),$ while $|\nabla_g \psi_j(x)|= \frac{2}{\rho_0}$ for
a.e. $x\in B_{x_j}(\rho_0)\setminus B_{x_j}(\frac{\rho_0}{2})$ (cf.
(\ref{dist-gradient})) and $|\nabla_g \psi_j(x)|= 0$ otherwise. The
uniform upper bound for the above covering yields that
\begin{equation}\label{no-ra-vonatkozo}
    1\leq \sum_{j\in \mathbb N}\psi_j(x)\leq N_0\ \ {\rm for\ all}\ x\in M.
\end{equation}

Let $\tau=\frac{4}{\rho_0}$ and fix $u\in C_0^\infty(M)$ arbitrarily
such that $\|u\|_{1,\tau}\leq 1.$ By the latter relation and the
properties of $\psi_j$ we have for every $j\in \mathbb N$ that
\begin{eqnarray*}
 \|\nabla_g(\psi_j^2u)\|_{L^n(M)}  &=& \|\psi_j^2\nabla_g
u+u\nabla_g\psi_j^2\|_{L^n(M)}\leq\|\psi_j^2\nabla_g
u\|_{L^n(M)}+2\|u\psi_j\nabla_g\psi_j\|_{L^n(M)} \\
   &\leq&\|\nabla_g
u\|_{L^n(M)}+\tau\|u\|_{L^n(M)}=\|u\|_{1,\tau}\\&\leq& 1.
\end{eqnarray*}
This estimate and relations (\ref{volume-unif-estimate}) and
(\ref{area-volume-333}) show that for every $j\in\mathbb N$ we can
apply Lemma \ref{lemma-isop}(i) to the geodesic ball
$B_{x_j}(\rho_0)$ and function $\psi_j^2u$ (standard density
arguments allow to consider that $\psi_j^2u$ is smooth), obtaining
for every $\alpha\in
\left[0,(C_0^nn^{-n}\omega_n^{-1})^\frac{1}{n-1}\alpha_n\right]$
that
\begin{equation}\label{becsles-elolrol}
    \int_{B_{x_j}(\rho_0)} \Phi_n(\alpha|\psi_j^2u|^\frac{n}{n-1}){\rm
d}v_g\leq M_0\eta_0\|\nabla_g
(\psi_j^2u)\|_{L^n(B_{x_j}(\rho_0))}^n.
\end{equation}
By the properties of the function $\psi_j$ and the covering of $M$,
it follows that
\begin{eqnarray*}
  \int_{M} \Phi_n(\alpha|u|^\frac{n}{n-1}){\rm d}v_g &\leq& \sum_{j\in \mathbb N}\int_{B_{x_j}(\frac{\rho_0}{2})} \Phi_n(\alpha|u|^\frac{n}{n-1}){\rm d}v_g \leq \sum_{j\in \mathbb N}\int_{B_{x_j}(\rho_0)} \Phi_n(\alpha|\psi_j^2u|^\frac{n}{n-1}){\rm d}v_g  \\
   &\leq&  M_0\eta_0\sum_{j\in \mathbb N}\|\nabla_g (\psi_j^2u)\|_{L^n(B_{x_j}(\rho_0))}^n\ \ \ \ \ \ \ \ \ \ \ \ \ \ \ \  \ \  \ \ \ \ \ \ \ \ \ \ \ \ \ \  \ \  \ \ \  {\rm [see\ (\ref{becsles-elolrol})]} \\
   &=& M_0\eta_0\sum_{j\in \mathbb N}\|\psi_j^2\nabla_g
u+u\nabla_g\psi_j^2\|_{L^n(B_{x_j}(\rho_0))}^n \\
   &\leq&M_0\eta_02^n\left(\sum_{j\in \mathbb N}\int_{B_{x_j}(\rho_0)}\psi_j|\nabla_g u|^n{\rm d}v_g+\frac{4^n}{\rho_0^n}\sum_{j\in \mathbb N}\int_{B_{x_j}(\rho_0)}\psi_j|u|^n{\rm
   d}v_g\right)\\&\leq& M_0\eta_02^nN_0\left(\int_{M}|\nabla_g u|^n{\rm d}v_g+\frac{4^n}{\rho_0^n}\int_{M}|u|^n{\rm
   d}v_g\right) \ \ \ \ \ \ \ \ \ \ \ \ \ \ \  \ \   {\rm [see\
   (\ref{no-ra-vonatkozo})]}\\&=&M_0\eta_02^nN_0\|u\|_{0,\tau}^n\\
   &\leq& M_0\eta_02^nN_0. \ \ \ \ \ \ \ \ \ \ \ \ \ \ \  \ \  \ \ \ \ \ \ \ \ \ \ \ \ \ \  \ \
   \ \ \ \ \ \ \ \
   \ [{\rm since}\ \|u\|_{0,\tau}\leq\|u\|_{1,\tau}\leq 1]
\end{eqnarray*}
Consequently, $S_{\alpha,\tau}^1(M,g)\leq M_0\eta_02^nN_0$ for
$\tau=\frac{4}{\rho_0}$ and every  $\alpha\in
\left[0,(C_0^nn^{-n}\omega_n^{-1})^\frac{1}{n-1}\alpha_n\right]$,
where the constants $M_0$, $C_0,$ $\eta_0,$ $N_0$ and $\rho_0$
depend only on $n,$ $k$ and $v.$

The continuity of the embedding $W^{1,n}(M)\hookrightarrow L^p(M)$
 for every $p\in [n,\infty)$ follows as in Yang \cite[Proposition 2.1]{Yang-JFA2012} whenever any of the assumptions (i), (ii) or (iii) holds.  \hfill
 $\square$

\begin{remark}\label{closer-and-closer}\rm If  ${\rm Isop}{(M,g)}$ is
close to 1, the constant $C_0$ in (\ref{area-volume-333}) can be
fixed close to $n\omega_n^\frac{1}{n}$. In such a case, the latter
proof shows that those numbers $\alpha\geq 0$ for which the
Moser-Trudinger inequality $({\bf MT})_{\alpha,\tau}^1$ holds on
$(M,g)$ approaches the critical exponent $\alpha_n$.
\end{remark}

{\it Proof of Corollary \ref{cor-Cartan Hadamard}.} If the
Cartan-Hadamard conjecture holds on $(M,g)$, then  ${\rm
 Isop}{(M,g)}=1.$
Now, if $\alpha\in
  \left[0,\alpha_n\right]$ and $\tau\geq
  1$, it remains to apply  Lemma \ref{lemma-isop} (ii).
The continuity of  embedding $W^{1,n}(M)\hookrightarrow L^p(M)$
 for every $p\in [n,\infty)$ follows again by  \cite[Proposition 2.1]{Yang-JFA2012}.
  \hfill
 $\square$

\begin{remark}\rm \label{becsles-Yang}
Let $(M,g)$ be an $n-$dimensional Hadamard manifold, $n\geq 2$.
Precisely as in the Euclidean case, one can prove:
\begin{equation}\label{alpha-minden}
    \Phi_n(\alpha|u|^\frac{n}{n-1})\in L^1(M)\ \ {\rm for\ all}\
\alpha>0,\ u\in W^{1,n}(M).
\end{equation}
The proof of (\ref{alpha-minden}) is based on the validity of the
Moser-Trudinger inequality $({\bf MT})_{\alpha,\tau}^0$ on $(M,g)$
for some $\alpha>0 $ and $\tau\geq
  1$  (cf.Corollary
\ref{cor-Cartan Hadamard}), the density of  $C_0^\infty(M)$ in
$W^{1,n}(M)$ endowed with the norm $\|\cdot\|_{0,\tau}$, and
 basic  properties of the function $\Phi_n$; a similar
 argument on Riemannian manifolds with Ricci curvature bounded from
 below is presented in Yang \cite[p. 1911]{Yang-JFA2012}.
\end{remark}

\vspace{0.5cm}

\section{An elliptic PDE with critical nonlinearity: proof of Theorem \ref{theorem-existence}}\label{sect-4}

Without mentioning explicitly, we assume throughout this section
that all assumptions of Theorem \ref{theorem-existence} are
satisfied. By $(f_0)$, one has that $f(0)=0$; therefore, we extend
continuously the function $f:[0,\infty)\to \mathbb R$ to the whole
$\mathbb R$ by $f(s)=0$ for $s\leq 0;$ thus,  $F(s)=0$ for $s\leq 0$
as well.  The function $u\in W^{1,n}(M)$ is a {\it weak solution} of
problem $(\mathcal P)$ if
\begin{equation}\label{weak-solution}
    \int_M (|\nabla_g u|^{n-2}\langle\nabla_g u, \nabla_g w\rangle+|u|^{n-2}uw){\rm
    d}v_g=\int_M f(u)w{\rm
    d}v_g\ \ {\rm for\ all}\ w\in W^{1,n}(M).
\end{equation}
By the above extension it turns out that every  weak solution of
problem $(\mathcal P)$ is non-negative.

 Let
$\mathcal E:W^{1,n}(M)\to \mathbb R$ be the energy functional
associated with problem $(\mathcal P)$, given by
$$\mathcal E(u)=\frac{\|u\|_{0,1}^n}{n}-\mathcal F(u),$$
where $$\mathcal F(u)=\int_M F(u){\rm d}v_g.$$  Due to $(f_0)$,
$(f_1)$, there exists $c_0>0$ such that
\begin{equation}\label{f-elso-becsles}
    |f(s)|\leq
c_0\left(|s|^{\gamma-1}+\Phi_n(\alpha_0|s|^\frac{n}{n-1})\right)\
{\rm for\ all}\ s\in \mathbb R.
\end{equation}
Therefore, by hypothesis $(f_2)$, H\"older's inequality and the
inequality
\begin{equation}\label{Phi-egyenlotlenseg}
    \Phi_n(s)^q\leq \Phi_n(qs)\ {\rm  for\ every}\ q\geq 1,
\end{equation}
 it
follows for every $u\in W^{1,n}(M)$ that
\begin{eqnarray}\label{F-re-vett-becsles}
  0\leq \mathcal F(u) &\leq& {c_0}\int_M|u|^{\gamma}{\rm d}v_g+c_0\int_M |u|\Phi_n(\alpha_0|u|^\frac{n}{n-1}){\rm d}v_g \nonumber\\
   &\leq& c_0\|u\|_{L^\gamma(M)}^\gamma+ c_0\|u\|_{L^n(M)}\left(\int_M \Phi_n\left(\frac{\alpha_0n}{n-1}|u|^\frac{n}{n-1}\right){\rm d}v_g\right)^\frac{n-1}{n}.\nonumber
\end{eqnarray}
The continuous embedding  $W^{1,n}(M)\hookrightarrow L^p(M)$ for
every $p\in [n,\infty)$ and relation (\ref{alpha-minden}) imply that
the latter term in the above estimate is finite, i.e., the energy
functional $\mathcal E$ is well-defined on $W^{1,n}(M);$
furthermore, $\mathcal E$ is of class $C^1$ on $W^{1,n}(M)$ and
standard arguments yield that the critical points of $\mathcal E$
are precisely the weak solutions of problem $(\mathcal P)$.

Let $G$ be a compact connected subgroup of ${\rm Isom}_g(M)$ with
the required properties, i.e., ${\rm Fix}_M(G)=\{x_0\}$ for some
$x_0\in M$ and ${\rm Card}(O^x_G)=\infty$ for every $x\in M\setminus
\{x_0\}$. The action of $G$ on $W^{1,n}(M)$ is defined by
\begin{equation}\label{action-of-the-group}
    (\sigma u)(x)=u(\sigma^{-1}(x)) \ \ {\rm for\ all}\ \sigma\in G,\ u\in
    W^{1,n}(M),\ x\in M,
\end{equation}
 where
$\sigma^{-1}:M\to M$ is the inverse of the isometry $\sigma$. Let
 $$W_G^{1,n}(M)=\{u\in W^{1,n}(M):\sigma u=u\
{\rm for\ all}\ \sigma\in G\}$$ be the subspace of $G-$invariant
functions of $W^{1,n}(M)$ and let $\mathcal E_G:W_G^{1,n}(M)\to
\mathbb R$ be the restriction of the energy functional $\mathcal E$
to $W_G^{1,n}(M)$.

Several lemmas are needed in order to complete the proof of Theorem
\ref{theorem-existence}.

\begin{lemma}\label{lemma-szimmetria-elv} Every critical point of  $\mathcal E_G$
is a non-negative $G-$invariant weak solution of  $(\mathcal P)$.
\end{lemma}

{\it Proof.} We first notice that  $G$ acts continuously on
$W^{1,n}(M)$ by relation (\ref{action-of-the-group}); for instance,
for every $\sigma_1,\sigma_2\in G$, $u\in W^{1,n}(M)$ and $x\in M$
one has
$$((\sigma_1\sigma_2)u)(x)=u((\sigma_1\sigma_2)^{-1}(x))=u(\sigma_2^{-1}(\sigma_1^{-1}(x)))=(\sigma_2u)(\sigma_1^{-1}(x))=(\sigma_1(\sigma_2u))(x),$$ while the other properties trivially hold.

We claim that $\mathcal E$ is $G-$invariant. To see this, let $u\in
W^{1,n}(M)$ and $\sigma\in G$ be arbitrarily fixed. Since
$\sigma:M\to M$ is an isometry on $M$, by
(\ref{action-of-the-group}), for every $x\in M$ we have
$$\nabla_g(\sigma u)(x)=D \sigma_{\sigma^{-1}(x)} \nabla_g u(\sigma^{-1}(x)),$$
where $D\sigma_{\sigma^{-1}(x)}:T_{\sigma^{-1}(x)}M\to T_x M$
denotes the differential of  $\sigma$ at the point $\sigma^{-1}(x)$.
Note that the (signed) Jacobian determinant of $\sigma$ is 1 and
$D\sigma_{\sigma^{-1}(x)}$ preserves inner products. Therefore, by
using the latter facts, relation (\ref{action-of-the-group}) and a
change of variables $y=\sigma^{-1}(x)$, it turns out that
\begin{eqnarray*}
  \|\sigma u\|_{0,1}^n &=& \int_M \left(|\nabla_g (\sigma u)(x)|_x^n + |(\sigma
u)(x)|^n\right){\rm d}v_g(x) \\
   &=&\int_M \left(|\nabla_g
u(\sigma^{-1}(x))|_{\sigma^{-1}(x)}^n +
|u(\sigma^{-1}(x))|^n\right){\rm d}v_g(x)=\int_M \left(|\nabla_g
u(y)|_y^n + |u(y)|^n\right){\rm d}v_g(y)\\&=&\|u\|_{0,1}^n,
\end{eqnarray*}
and
$$\mathcal F(\sigma u)= \int_M F((\sigma u)(x)){\rm d}v_g(x)=\int_M F(u(\sigma^{-1} (x))){\rm d}v_g(x)=\int_M F(u(y)){\rm d}v_g(y)=\mathcal F(u),$$
which ends the proof of the claim.

Note that ${\rm Fix}_{W^{1,n}(M)}(G)$ is nothing but $W_G^{1,n}(M)$;
therefore,  if $u_G\in W_G^{1,n}(M)$ is a critical point of
$\mathcal E_G$, then due to Proposition \ref{Palais-PSC}, $u_G$ is
also a critical point of $\mathcal E$ and as such,  $u_G$ turns out
to be a $G-$invariant non-negative weak solution of $(\mathcal P)$,
as we pointed out before.  \hfill
 $\square$

\begin{lemma}\label{mountain-pass-lemma} The functional $\mathcal
E_G$ has the mountain pass geometry, i.e.,
\begin{itemize}
  \item[{\rm (i)}] for every non-negative, compactly supported $\tilde u\in W_G^{1,n}(M)\setminus \{0\}$ we have $\mathcal
E_G(s\tilde u)\to -\infty$ as $s\to \infty;$
  \item[{\rm (ii)}] there exist $\tilde r>0$ and $\tilde \delta>0$ such that $\mathcal
E_G(u)\geq\tilde  \delta$ for every $u\in
  W_G^{1,n}(M)$ with $\|u\|_{0,1}=\tilde r$.
\end{itemize}
\end{lemma}

{\it Proof.} (i) Let $\tilde u\in W_G^{1,n}(M)\setminus\{0\}$ be a
non-negative function with compact support contained in the geodesic
ball $\overline B_{x_0}(r)$ for some $r>0$. By $(f_2),$ it follows
that there exist $c_1,c_2>0$ such that $F(t)\geq c_1 t^\mu-c_2$ for
every $t\in [0,\infty).$ Therefore,
$$\mathcal E_G(s\tilde u)=s^n\frac{\| \tilde u\|_{0,1}^n}{n}-\mathcal F(s\tilde u)\leq
s^n\frac{\| \tilde u\|_{0,1}^n}{n}-c_1s^\mu\int_{\overline
B_{x_0}(r)}  \tilde u^\mu{\rm d}v_g +c_2{\rm Vol}_g(\overline
B_{x_0}(r)).$$ Since $\tilde u\neq 0$ and $\mu>n$, one has that
$\mathcal E_G(s \tilde u)\to -\infty$ as $s\to \infty.$

(ii) By $(f_0)$ and $(f_1)$, there exists  $c_3>0$ such that
\begin{equation}\label{majd-kesobbi-f-becsles}
    |f(s)|\leq
c_3|s|^{\gamma-1}\left(1+\Phi_n(\alpha_0|s|^\frac{n}{n-1})\right)\
{\rm for\ all}\ s\in \mathbb R.
\end{equation}
 By H\"older's inequality and (\ref{Phi-egyenlotlenseg}), for every
$u\in W_G^{1,n}(M)$  one has
\begin{eqnarray}\label{F-rol-kell}
  \mathcal F(u) &\leq& {c_3}\|u\|^{\gamma}_{L^\gamma(M)}+c_3\int_M
|u|^\gamma\Phi_n(\alpha_0|u|^\frac{n}{n-1}){\rm d}v_g \nonumber\\
 &\leq&{c_3}\|u\|^{\gamma}_{L^\gamma(M)}+c_3\|u\|^{\gamma}_{L^{2\gamma}(M)}\left(\int_M
\Phi_n(2\alpha_0|u|^\frac{n}{n-1}){\rm d}v_g\right)^\frac{1}{2}.
\end{eqnarray}
Due to \cite[Theorem 1.2]{YSK} (or Corollary \ref{cor-Cartan Hadamard}
in dimensions 2, 3 and 4), the Moser-Trudinger inequality $({\bf
MT})_{\alpha_n,1}^0$ is valid on $(M,g)$, i.e.,
$S_{\alpha_n,1}^0(M,g)<\infty.$ Let $\mathfrak{s}_p>0$ be the best
embedding constant in $W^{1,n}(M)\hookrightarrow L^p(M)$, $p\in
[n,\infty)$, and let us choose $\tilde r>0$ such that
\begin{equation}\label{muszaly}
    2\alpha_0 \tilde r^\frac{n}{n-1}\leq\alpha_n\ {\rm and}\
c_3n\left(\mathfrak{s}_{\gamma}^\gamma+\mathfrak{s}_{2\gamma}^\gamma
(S_{ \alpha_n,1}^0(M,g))^\frac{1}{2}\right)\tilde r^{\gamma-n}<1.
\end{equation}
 Thus, for every $u\in W_G^{1,n}(M)$ with
$\|u\|_{0,1}=\tilde r$, by relations (\ref{F-rol-kell}) and
(\ref{muszaly}) it follows that
$$\mathcal E_G(u)\geq \frac{\tilde r^n}{n}-c_3\left(\mathfrak{s}_{\gamma}^\gamma+\mathfrak{s}_{2\gamma}^\gamma (S_{\alpha_n,1}^0(M,g))^\frac{1}{2}\right)\tilde r^\gamma:=\tilde \delta>0,$$
which concludes the proof. \hfill
 $\square$\\

 The next lemma gives information on the behavior of
 Palais-Smale sequences of the functional $\mathcal
E_G$; let $W_G^{1,n}(M)^*$ be the dual of $W_G^{1,n}(M),$ and
$\langle\cdot, \cdot\rangle_*$ be the duality pairing between
$W_G^{1,n}(M)^*$ and $W_G^{1,n}(M).$

\begin{lemma}\label{Palais-Smale} If $\{u_j\}_{j\in \mathbb N}\subset W_G^{1,n}(M)$ is a Palais-Smale sequence of $\mathcal
E_G$, i.e., $\mathcal E_G(u_j)\to c\in \mathbb R$ and $\mathcal
E_G^\prime(u_j)\to 0$ in $W_G^{1,n}(M)^*,$ then there exist a
subsequence of $\{u_j\}$ $($still denoted by $\{u_j\})$ and $u_G\in
W_G^{1,n}(M)$ such that
\begin{itemize}
  \item[{\rm (i)}] $\lim_{j\to \infty}\mathcal F(u_j)= \mathcal F(u_G);$
  \item[{\rm (ii)}]  $u_j\to u_G\ {\rm strongly\ in}\ L^p(M)\ {\rm for\ every}\ p\in
(n,\infty);$
  \item[{\rm (iii)}] $\mathcal E_G^\prime(u_G)=0,$ i.e., $u_G$ is a
  critical point of $\mathcal E_G$.
\end{itemize}
\end{lemma}

{\it Proof.} (i)\&(ii) Let $\{u_j\}_{j\in \mathbb N}\subset
W_G^{1,n}(M)$ be a Palais-Smale sequence of $\mathcal E_G$ at level
$c\in \mathbb R$, i.e., $\mathcal E_G(u_j)\to c$ and $|\langle
\mathcal E_G^\prime(u_j),w \rangle_*| \leq \varepsilon_j\|w\|_{0,1}$
for every $w\in
    W_G^{1,n}(M)$, where $\lim_{j\to \infty}\varepsilon_j= 0;$
    explicitly, one has
\begin{equation}\label{PS-1}
\frac{\|u_j\|_{0,1}^n}{n}-\mathcal F(u_j)\to c;
\end{equation}
\begin{equation}\label{PS-2}
\left|\int_M (|\nabla_g u_j|^{n-2}\langle\nabla_g u_j,\nabla_g
w\rangle+|u_j|^{n-2}u_jw){\rm
    d}v_g-\int_M f(u_j)w{\rm
    d}v_g\right|\leq \varepsilon_j\|w\|_{0,1}, {\rm \forall} w\in
    W_G^{1,n}(M).
\end{equation}
By construction,  $f(s)=F(s)=0$ for $s\leq 0$; thus, multiplying
relation (\ref{PS-1}) by $\mu$, letting $w=u_j$ in (\ref{PS-2}), and
adding these relations, it follows by hypothesis $(f_2)$  that
$$\left(\frac{\mu}{n}-1\right)\|u_j\|_{0,1}^n\leq\int_M (\mu F(u_j)-f(u_j)u_j){\rm
    d}v_g+\mu|c|+\varepsilon_j\|u_j\|_{0,1}\leq \mu|c|+\varepsilon_j\|u_j\|_{0,1}.$$
Since $\mu>n$, the sequence $\{u_j\}$ is bounded in $W_G^{1,n}(M)$;
in particular, by relation (\ref{PS-1}) and the latter estimate one
can guarantee the existence of  $c_4>0$ (depending only on $n,$
$\mu$ and $c)$ such that for every $j\in \mathbb N$,
\begin{equation}\label{f-ek-becslese}
   \mathcal F(u_j)= \int_M F(u_j){\rm
    d}v_g\leq c_4\ \ {\rm and}\ \ \int_M f(u_j)u_j{\rm
    d}v_g\leq c_4.
\end{equation}
By the boundedness of $\{u_j\}$ in $W_G^{1,n}(M)$ together with the
hypothesis Fix$_M(G)=\{x_0\}$ and Proposition \ref{Tintarev}, there
exists $u_G\in W_G^{1,n}(M)$ such that, up to a subsequence, we have
\begin{equation}\label{u_j-weak}
     u_j\rightharpoonup u_G\ {\rm weakly\ in}\  W_G^{1,n}(M);
\end{equation}
\begin{equation}\label{u_j-strong}
    u_j\to u_G\ {\rm strongly\ in}\ L^p(M)\ {\rm for\ every}\ p\in
(n,\infty);
\end{equation}
\begin{equation}\label{u_j-a.e.}
     u_j\to u_G\ {\rm a.e.\ in}\  M.
\end{equation}
Let $\varepsilon>0$ be fixed arbitrarily, and let
\begin{equation}\label{K-ra_becsles}
    K>\max\left\{R_0,\frac{A_0}{\varepsilon}
c_4,\frac{A_0}{\varepsilon}\int_M f(u_G)u_G{\rm
    d}v_g\right\},
\end{equation}
where $R_0>0$ and $A_0>0$ are from $(f_3)$. Since $F(s)=0$ for every
$s\in (-\infty,0]$ and $f(s)s\geq 0$ for every $s\in [0,\infty)$
(cf. $(f_2)$), by hypothesis $(f_3)$ and relations
(\ref{K-ra_becsles}) and (\ref{f-ek-becslese}), one has for every
$j\in \mathbb N$ that
\begin{eqnarray}\label{F-eps-becsles-1}
\nonumber  \int_{\{|u_j|>K\}} F(u_j){\rm d}v_g &=& \int_{\{u_j>K\}}
F(u_j){\rm d}v_g\leq A_0\int_{\{u_j>K\}} f(u_j){\rm d}v_g\\&\leq&
\frac{A_0}{K}\int_{\{u_j>K\}} f(u_j)u_j{\rm
d}v_g\leq \frac{A_0}{K}c_4\nonumber \\
   &<&\varepsilon.
\end{eqnarray}
In a similar way, we have
\begin{eqnarray}\label{F-eps-becsles-2}
  \int_{\{|u_G|>K\}} F(u_G){\rm d}v_g \leq
A_0\int_{\{u_G>K\}} f(u_G){\rm d}v_g\leq
\frac{A_0}{K}\int_{\{u_G>K\}} f(u_G)u_G{\rm d}v_g<\varepsilon.
\end{eqnarray}
By relation (\ref{majd-kesobbi-f-becsles}),  it follows that
$f(s)\leq
c_3s^{\gamma-1}\left(1+\Phi_n(\alpha_0K^\frac{n}{n-1})\right)\ {\rm
for\ all}\ s\in [0,K].$ Therefore,
$$F(s)\leq
c_5s^{\gamma}\ {\rm for\ all}\ s\in [0,K],$$ where
$c_5=c_3\left(1+\Phi_n(\alpha_0K^\frac{n}{n-1})\right).$
Consequently, for every $j\in \mathbb N$ we have
\begin{equation}\label{Lebesgue-uj}
    \chi_{\{|u_j|\leq K\}}F(u_j)\leq c_5 |u_j|^\gamma,
\end{equation}
where $\chi_A$ denotes the characteristic function of the set
$A\subset M.$ We recall   the inequality
\begin{equation}\label{egyenlotlenseg}
    \left||s|^p-|t|^p\right|\leq p|s-t|(|s|^{p-1}+|t|^{p-1})\ {\rm for\
    all}\ p>1\ {\rm and}\ t,s\in \mathbb R.
\end{equation}
 By (\ref{egyenlotlenseg}) and H\"older's inequality, one
has
\begin{eqnarray*}
  \int_M\left||u_j|^\gamma-|u_G|^\gamma\right|{\rm d}v_g &\leq& \gamma\int_M|u_j-u_G|(|u_j|^{\gamma-1}+|u_G|^{\gamma-1}){\rm
d}v_g \\
   &\leq&\gamma\|u_j-u_G\|_{L^\gamma(M)}(\|u_j\|^{\gamma-1}_{L^\gamma(M)}+\|u_G\|^{\gamma-1}_{L^\gamma(M)}).
\end{eqnarray*}
Since $\gamma>n$, due to (\ref{u_j-strong}) the latter term tends to
zero, thus $|u_j|^\gamma$ converges to $|u_G|^\gamma$ in $L^1(M)$ as
$j\to \infty.$ By (\ref{u_j-a.e.}), (\ref{Lebesgue-uj}) and the
generalized Lebesgue dominated convergence theorem we have
$$\lim_{j\to \infty}\int_M\chi_{\{|u_j|\leq K\}}F(u_j){\rm d}v_g=\int_M\chi_{\{|u_G|\leq K\}}F(u_G){\rm d}v_g.$$
The latter relation together with (\ref{F-eps-becsles-1}) and
(\ref{F-eps-becsles-2}) implies that
$$\lim_{j\to \infty}\int_MF(u_j){\rm d}v_g=\int_MF(u_G){\rm d}v_g,$$
which proves (i). Note that (\ref{u_j-strong}) is precisely the
property (ii).

 (iii) The proof is divided into
several steps.

 \underline{\it Step 1}:
\begin{equation}\label{fuu}
\lim_{j\to \infty}\int_Mf(u_j)w{\rm d}v_g=\int_Mf(u_G)w{\rm d}v_g\
{\rm for\ all}\ w\in C_0^\infty(M).
\end{equation}
This step is similar to (i); let $\varepsilon>0$ and $w\in
C_0^\infty(M)\setminus \{0\}$ be arbitrarily fixed, and let
$$
    K>\frac{\|w\|_{L^\infty(M)}}{\varepsilon}\max\left\{
c_4,\int_M f(u_G)u_G{\rm
    d}v_g\right\}.
$$
Relation (\ref{f-ek-becslese}), the choice of $K>0$ and the fact
that $|f(s)s|=f(s)s$ for every $s\in \mathbb R$ show that
\begin{equation}\label{2-becsles-rend}
    \int_{\{|u_j|>K\}} \left|f(u_j)w\right|{\rm d}v_g<\varepsilon\ \
{\rm and}\ \ \int_{\{|u_G|>K\}} \left|f(u_G)w\right|{\rm
d}v_g<\varepsilon.
\end{equation}
As above, by (\ref{majd-kesobbi-f-becsles}),  one has $f(s)\leq
c_3s^{\gamma-1}\left(1+\Phi_n(\alpha_0K^\frac{n}{n-1})\right)\ {\rm
for\ all}\ s\in [0,K].$ Therefore,
\begin{equation}\label{Lebesgue-uj-2}
    \chi_{\{|u_j|\leq K\}}|f(u_j)w|\leq c_6 |u_j|^{\gamma-1}|w|,
\end{equation}
where $c_6=c_3\left(1+\Phi_n(\alpha_0K^\frac{n}{n-1})\right),$ which
is formally the same as $c_5$ but perhaps $K$ differs. Note that
$|u_j|^{\gamma-1}|w|$ converges to $|u_G|^{\gamma-1}|w|$ in
$L^1(M)$; indeed, since $\gamma>n\geq 2,$ by (\ref{egyenlotlenseg})
and H\"older's inequality we have
\begin{eqnarray*}
  \int_M\left||u_j|^{\gamma-1}-|u_G|^{\gamma-1}\right||w|{\rm d}v_g &\leq& (\gamma-1)\int_M|u_j-u_G|(|u_j|^{\gamma-2}+|u_G|^{\gamma-2})|w|{\rm
d}v_g \\
   &\leq&(\gamma-1)\|u_j-u_G\|_{L^\gamma(M)}(\|u_j\|^{\gamma-2}_{L^\gamma(M)}+\|u_G\|^{\gamma-2}_{L^\gamma(M)})\|w\|_{L^\gamma(M)},
\end{eqnarray*}
and according to (\ref{u_j-strong}), the above integral tends to
zero as $j\to \infty.$ The generalized Lebesgue dominated
convergence theorem together with (\ref{u_j-a.e.}) and
(\ref{Lebesgue-uj-2}) provide
$$\lim_{j\to \infty}\int_M\chi_{\{|u_j|\leq K\}}f(u_j)w{\rm d}v_g=\int_M\chi_{\{|u_G|\leq K\}}f(u_G)w{\rm d}v_g.$$
Combining the latter relation with (\ref{2-becsles-rend}), the claim
(\ref{fuu}) follows.

 \underline{\it Step 2}: {\it for every compact set $S\subset M\setminus \{x_0\}$,
 one has}
 \begin{equation}\label{fuu-2}
\lim_{j\to \infty}\int_S|f(u_j)(u_j-u_G)|{\rm d}v_g=0.
\end{equation}
In order to prove this claim, let $\delta_0>0$ be fixed such that
\begin{equation}\label{delta-null-valasztasa}
    \alpha_0 \frac{\gamma}{\gamma-1}2^\frac{n+1}{n-1}\delta_0^\frac{1}{n-1}<\alpha_n,
\end{equation}
where $\gamma>n$ and $\alpha_0>0$  are from hypotheses $(f_0)$ and
$(f_1),$ respectively. We are going to prove first an
energy-concentration property; namely, we claim that for every $x\in
M\setminus \{x_0\}$ there exists $0<r_x<d_g(x_0,x)$ such that
\begin{equation}\label{energia-koncentracio}
    \lim_{j\to \infty}\int_{B_x(r_x)}(|\nabla_g u_j|^n+|u_j|^n){\rm
d}v_g<\delta_0.
\end{equation}
By contradiction, we assume that there exists $\tilde x\in
M\setminus \{x_0\}$ such that $$\lim_{r\to 0}\lim_{j\to
\infty}\int_{B_{\tilde x}(r)}(|\nabla_g u_j|^n+|u_j|^n){\rm
d}v_g\geq \delta_0.$$ By assumption, we have ${\rm Card}(O^{\tilde
x}_G)=\infty$; thus, we may fix the distinct points $\tilde
x_1,...,\tilde x_N\in O^{\tilde x}_G$  with
$$N>\frac{n(|c|+c_4)}{\delta_0},$$ where $c\in \mathbb R$ and
$c_4>0$ are from (\ref{PS-1}) and (\ref{f-ek-becslese}),
respectively. Note that there exists $\sigma_l\in G$ such that
$\tilde x_l=\sigma_l( \tilde x)$ for every $l\in \{1,...,N\}.$
Furthermore,  $B_{\tilde x_l}(r)=\sigma_l B_{\tilde x}(r)$ for every
$l\in \{1,...,N\}.$ By using these facts, since $u_j$ are
$G-$invariant functions and $\sigma_l\in G$ are isometries on $M$, a
similar argument as in the proof of Lemma \ref{lemma-szimmetria-elv}
shows that for every $l\in \{1,...,N\},$
$$\int_{B_{\tilde x_l}(r)}(|\nabla_g u_j|^n+|u_j|^n){\rm
d}v_g=\int_{\sigma_lB_{\tilde x}(r)}(|\nabla_g u_j|^n+|u_j|^n){\rm
d}v_g=\int_{B_{\tilde x}(r)}(|\nabla_g u_j|^n+|u_j|^n){\rm d}v_g.$$
By relations (\ref{PS-1}), (\ref{f-ek-becslese}) and the above
assumption, it follows that
\begin{eqnarray*}
  {n(|c|+c_4)} &\geq& \lim_{j\to \infty}\|u_j\|_{0,1}^n=\lim_{j\to
\infty}\int_M(|\nabla_g u_j|^n+|u_j|^n){\rm d}v_g \\
   &\geq&\sum_{l=1}^N\lim_{r\to 0}\lim_{j\to \infty}\int_{B_{\tilde
x_l}(r)}(|\nabla_g u_j|^n+|u_j|^n){\rm d}v_g=N\lim_{r\to
0}\lim_{j\to
\infty}\int_{B_{\tilde x}(r)}(|\nabla_g u_j|^n+|u_j|^n){\rm d}v_g\\
&\geq& N\delta_0,
\end{eqnarray*}
which contradicts the choice of $N.$ Therefore, relation
(\ref{energia-koncentracio}) holds.

Let $x\in M\setminus \{x_0\}$ be arbitrarily fixed, $r:=r_x>0$ from
(\ref{energia-koncentracio}) and  $a_j=\frac{1}{{\rm
Vol}_g(B_{x}(r))}\int_{B_{x}(r)} u_j {\rm d}v_g$. By H\"older's
inequality and (\ref{energia-koncentracio}), for enough large $j\in
\mathbb N$ we have
$$|a_j|\leq \frac{1}{{\rm
Vol}_g(B_{x}(r))}\int_{B_{x}(r)} |u_j| {\rm d}v_g\leq {{\rm
Vol}_g(B_{x}(r))}^{-\frac{1}{n}}\left(\int_{B_{x}(r)} |u_j|^n {\rm
d}v_g\right)^\frac{1}{n}\leq\left(\frac{\delta_0}{{{\rm
Vol}_g(B_{x}(r))}}\right)^\frac{1}{n}.$$ Let $\tilde u_j=u_j-a_j$
for every $j\in \mathbb N.$ Then for enough large $j\in \mathbb N$,
one has
$$\int_{B_x(r)}\tilde u_j{\rm d}v_g=0\ \ {\rm and}\ \
\int_{B_x(r)}|\nabla_g\tilde u_j|^n{\rm d}v_g<\delta_0.$$ Therefore,
by relation (\ref{delta-null-valasztasa}) and Cherrier's result (cf.
(\ref{Cherrier-egyenlotlenseg})) applied on $\overline B_x(r)$ for
the functions $\frac{\tilde u_j}{\|\nabla_g \tilde
u_j\|_{L^n(B_x(r))}}$, $j\in \mathbb N$ large enough, it follows
that
\begin{eqnarray*}
 \int_{B_x(r)}e^{\alpha_0\frac{\gamma}{\gamma-1}|u_j|^\frac{n}{n-1}}{\rm d}v_g  &=& \int_{B_x(r)}e^{\alpha_0\frac{\gamma}{\gamma-1}|\tilde u_j+a_j|^\frac{n}{n-1}}{\rm d}v_g
   \leq e^{\alpha_0\frac{\gamma}{\gamma-1}2^\frac{n}{n-1}|a_j|^\frac{n}{n-1}}\int_{B_x(r)}e^{\alpha_0\frac{\gamma}{\gamma-1}2^\frac{n}{n-1}|\tilde u_j|^\frac{n}{n-1}}{\rm d}v_g  \\
   &\leq &c_7,
\end{eqnarray*}
where the constant $c_7>0$ depends on $\alpha_0,$ $n,$ $r,$ $x,$
$\gamma$ and $\delta_0$, but not on $j\in \mathbb N.$

Since $\{u_j\}$ is bounded in $L^\gamma(M)$, the latter estimate
together with  H\"older's inequality and relations
(\ref{f-elso-becsles}) and (\ref{Phi-egyenlotlenseg}) yield
\begin{eqnarray*}
  I_j&:=&\int_{B_x(r)}|f(u_j)(u_j-u_G)|{\rm d}v_g \leq \left(\int_{B_x(r)}|f(u_j)|^\frac{\gamma}{\gamma-1}{\rm d}v_g\right)^{1-\frac{1}{\gamma}}\left( \int_{B_x(r)}|u_j-u_G|^\gamma{\rm d}v_g\right)^\frac{1}{\gamma} \\
   &\leq& 2c_0\left(\int_{B_x(r)}|u_j|^{\gamma}{\rm d}v_g+\int_{B_x(r)}\Phi_n\left(\alpha_0\frac{\gamma}{\gamma-1}|u_j|^\frac{n}{n-1}\right){\rm d}v_g\right)^{1-\frac{1}{\gamma}}\|u_j-u_G\|_{L^\gamma({B_x(r)})} \\
   &\leq&c_8\|u_j-u_G\|_{L^\gamma(M)},
\end{eqnarray*}
where  $c_8>0$ does not depend on $j\in \mathbb N.$ Consequently,
due to (\ref{u_j-strong}), we have $$\lim_{j\to \infty} I_j=0.$$
Now, the compact set  $S\subset M\setminus \{0\}$ can be covered by
a finite number of geodesic balls with the above properties, which
completes the proof of (\ref{fuu-2}) throughout the latter limit.

 \underline{\it Step 3}: {\it for every compact set $S\subset M\setminus \{x_0\}$,
 one has}
 \begin{equation}\label{fuu-3}
\lim_{j\to
\infty}\int_S\left(|\nabla_gu_j-\nabla_gu_G|^n+|u_j-u_G|^n\right){\rm
d}v_g=0.
\end{equation}
Let $x\in M\setminus \{x_0\}$ be arbitrarily fixed  and
$r:=r_x<d_g(x_0,x)$ from (\ref{energia-koncentracio}). For every
$0<\rho\leq r$, let $A_{x_0}(\rho)=B_{x_0}(d_g(x_0,x)+\rho)\setminus
\overline B_{x_0}(d_g(x_0,x)-\rho)$ be the open geodesic annulus
with center $x_0\in M$ and radii $d_g(x_0,x)\pm\rho$, respectively.

We consider a $d_g(x_0,\cdot)-$radially symmetric function
$\varphi\in C_0^\infty(A_{x_0}(r))$ such that $0\leq \varphi\leq 1$
and $\varphi=1$ on $A_{x_0}(\frac{r}{2}).$ Hereafter, a function
$\varphi:M\to \mathbb R$ is called $d_g(x_0,\cdot)-$radially
symmetric, if there exists a function $h_\varphi:[0,\infty)\to
\mathbb R$ such that $\varphi(x)=h_\varphi(d_g(x_0,x))$ for every
$x\in M.$ For simplicity, we extend $\varphi$ by zero to the whole
$M$ outside of the geodesic annulus $A_{x_0}(r)$.

Note that $\varphi$ is $G-$invariant. Indeed, since ${\rm
Fix}_M(G)=\{x_0\}$, for every $x\in M$ and isometry $\sigma\in G$ we
have
$$\varphi(\sigma (x))=h_\varphi(d_g(x_0,\sigma(x)))=h_\varphi(d_g(\sigma(x_0),\sigma(x)))=h_\varphi(d_g(x_0,x))=\varphi(x).$$
In particular,  $\varphi(u_j-u_G)\in W_G^{1,n}(M)$ for every $j\in
\mathbb N$; insert this test-function into (\ref{PS-2}) to obtain
$$\int_M |\nabla_g u_j|^{n-2}\langle\nabla_g u_j,(u_j-u_G)\nabla_g\varphi+\varphi (\nabla_g u_j-\nabla_g u_G)\rangle{\rm
    d}v_g
+\int_M \varphi|u_j|^{n-2}u_j(u_j-u_G){\rm
    d}v_g$$ $$-\int_M \varphi f(u_j)(u_j-u_G){\rm
    d}v_g\leq \varepsilon_j\|\varphi(u_j-u_G)\|_{0,1}.$$
 Reorganizing this inequality, it yields that
\begin{eqnarray*}
  J_j &:=& \int_{A_{x_0}(r)} \varphi\langle|\nabla_g u_j|^{n-2}\nabla_g u_j-|\nabla_g u_G|^{n-2}\nabla_g u_G, \nabla_g u_j-\nabla_g u_G\rangle{\rm
    d}v_g
\\&&+\int_{A_{x_0}(r)} \varphi\left(|u_j|^{n-2}u_j-|u_G|^{n-2}u_G\right)(u_j-u_G){\rm
    d}v_g \\
    &\leq& \int_{A_{x_0}(r)}(u_G-u_j)|\nabla_g u_j|^{n-2}\langle\nabla_g u_j,\nabla_g\varphi\rangle{\rm d}v_g+\int_{A_{x_0}(r)}\varphi|\nabla_g u_G|^{n-2}\langle\nabla_g u_G, \nabla_g u_G-\nabla_g u_j\rangle{\rm
    d}v_g\\
    &&+\int_{A_{x_0}(r)} \varphi|u_G|^{n-2}u_G(u_G-u_j){\rm
    d}v_g+\int_{A_{x_0}(r)} \varphi f(u_j)(u_j-u_G){\rm
    d}v_g+ \varepsilon_j\|\varphi(u_j-u_G)\|_{0,1}.
\end{eqnarray*}
We shall check that every term on the right hand side of the above
inequality tend to $0$ as $j\to \infty$. First, by H\"older's
inequality, we have $$ \left|\int_{A_{x_0}(r)}(u_j-u_G)|\nabla_g
u_j|^{n-2}\langle\nabla_g u_j,\nabla_g\varphi\rangle{\rm
d}v_g\right|\leq \ \ \ \ \ \ \ \ \ \ \ \ \ \  \ \ \ \ \ \ \ \ \ \ \
\ \ \ \ \ \ \ \ \ \ \ \ \ \ \ \ \ \ \
   $$$$ \ \ \ \ \ \ \ \ \ \ \ \ \ \ \ \ \ \ \ \ \ \ \ \ \ \ \ \ \leq \|u_j-u_G\|_{L^\gamma(M)}{\rm Vol}_g(A_{x_0}(r))^\frac{\gamma-n}{\gamma
   n}\|\nabla_g
   u_j\|^{n-1}_{L^n(M)}\|\nabla_g\varphi\|_{L^\infty(M)}.
$$
Since $\{u_j\}$ is bounded in $W_G^{1,n}(M)$ and $\gamma>n$, due to
(\ref{u_j-strong}), the latter expression tends to $0$ as  $j\to
\infty$.
Second, due to (\ref{u_j-weak}), one has in particular that
$\nabla_g u_j\rightharpoonup \nabla_g u_G$ weakly in
$L^n(A_{x_0}(r),TM)$.
 Therefore, $$\lim_{j\to
\infty}\int_{A_{x_0}(r)}\langle\varphi|\nabla_g u_G|^{n-2}\nabla_g
u_G,\nabla_g u_j-\nabla_g u_G\rangle{\rm
    d}v_g=0.$$
The third term trivially converges to $0.$ Due to (\ref{fuu-2}), the
fourth term tends to $0$ as well. Since $\{\varphi(u_j-u_G)\}$ is
bounded in $W_G^{1,n}(M)$ and $\lim_{j\to \infty}\varepsilon_j=0$,
the latter term on the right hand side also tends to $0$.
Consequently,
\begin{equation}\label{J-j-becsles}
    \lim_{j\to \infty}J_j\leq 0.
\end{equation}
 On the other hand, for every $x\in M$ and $X,Y\in T_x M$, we have the
inequality
$$2^{2-n}|X-Y|^n \leq \langle |X|^{n-2}X-|Y|^{n-2}Y,X-Y\rangle.$$
Combining this inequality with (\ref{J-j-becsles}) and using the
properties of $\varphi$, it turns out that
$$ \lim_{j\to \infty}\int_{A_{x_0}(\frac{r}{2})}\left(|\nabla_g u_j-\nabla_g u_G|^n+|u_j-u_G|^n\right){\rm
    d}v_g=0.$$
It remains to apply a covering argument as in Step 2 in order to
prove (\ref{fuu-3}).

\underline{\it Step 4}: {\it  concluding the proof}. By Step 3 (cf.
(\ref{fuu-3})), we get in particular that the sequence $\{\nabla_g
u_j\}$ converges (up to a subsequence) to $\nabla_g u_G$ almost
everywhere on $M.$ Since the sequence $\{|\nabla_g
u_j|^{n-2}\nabla_g u_j\}$ is bounded in $L^\frac{n}{n-1}(M,TM)$,
there exists $X_0\in TM$ such that $|\nabla_g u_j|^{n-2}\nabla_g
u_j\rightharpoonup X_0$ weakly in $L^\frac{n}{n-1}(M,TM)$. The a.e.
convergence of the sequence $\{\nabla_g u_j\}$ to $\nabla_g u_G$
implies that $X_0$ should be precisely $|\nabla_g u_G|^{n-2}\nabla_g
u_G$. Consequently,
\begin{equation}\label{nabla-weak}
    |\nabla_g u_j|^{n-2}\nabla_g u_j\rightharpoonup |\nabla_g
u_G|^{n-2}\nabla_g u_G\ {\rm weakly\ in}\ L^\frac{n}{n-1}(M,TM).
\end{equation}

Let $w\in W_G^{1,n}(M)$ be arbitrarily fixed. By density, there
exists a sequence $\{w_l\}\subset C_0^\infty(M)$ which converges to
$w$ in $\|\cdot\|_{0,1}$. By using $w_l$ as a test-function in
(\ref{PS-2}), due to relations (\ref{fuu}), (\ref{nabla-weak}) and
the fact that $\lim_{j\to \infty}\varepsilon_j=0,$ we have
$$\int_M (|\nabla_g u_G|^{n-2}\langle\nabla_g u_G,\nabla_g
w_l\rangle+|u_G|^{n-2}u_Gw_l){\rm
    d}v_g-\int_M f(u_G)w_l{\rm
    d}v_g=0\ {\rm for\ all}\ l\in \mathbb N.$$
Letting now $l\to \infty$, it turns out that
$$\int_M (|\nabla_g u_G|^{n-2}\langle\nabla_g u_G,\nabla_g
w\rangle+|u_G|^{n-2}u_Gw){\rm
    d}v_g-\int_M f(u_G)w{\rm
    d}v_g=0,$$
which is nothing but $\langle\mathcal E_G^\prime(u_G),w\rangle_*=0;$
thus, the arbitrariness of $w\in W_G^{1,n}(M)$ implies that
$\mathcal E_G^\prime(u_G)=0,$ concluding the proof.
 \hfill
 $\square$\\

 Since $(M,g)$ is a Hadamard manifold,  its injectivity radius is $+\infty$; thus, it costs no generality to consider  in particular $\varepsilon_0=1$ and
$\varepsilon:=\frac{1}{j}$ $(j\in \mathbb N\setminus \{1\})$ in the
function (\ref{Moser-function-1}), introducing the rescaled Moser
functions
\begin{equation}\label{Moser-function-2}
   m_j(x):=\frac{(\log
    j)^\frac{n-1}{n}}{\omega_{n-1}^\frac{1}{n}}u_\frac{1}{j}(x)=\frac{(\log
    j)^\frac{n-1}{n}}{\omega_{n-1}^\frac{1}{n}}\min\left\{\left(-\frac{\log{d_g(x_0,x)}}{\log
    j}\right)_+,1\right\},\ x\in M.
\end{equation}
 The functions $m_j$
are well-defined and supp$(m_j)=\overline B_{x_0}(1)$ for every
$j\in \mathbb
 N\setminus\{1\}.$ Moreover, since ${\rm Fix}_M(G)=\{x_0\}$, it follows that the functions $m_j$ are $G-$invariant for every $j\in \mathbb N\setminus
\{1\}$; thus, $m_j\in W_G^{1,n}(M)$. Taking into account the
computations from the proof of Proposition
\ref{prop-veges-terfogat}, it follows that
\begin{equation}\label{m_j-aszimp}
    \|m_j\|_{0,1}^n=1+\mathcal O\left(\frac{1}{\log j}\right)\ {\rm as}\
j\to \infty.
\end{equation}

Moreover, inspired by Adimurthi and  Yang \cite{Adimurthi-Yang} and
do \'O \cite{doO}, we have

\begin{lemma}\label{Moser-fuggvenyek-fo} There exists $j_0\in \mathbb N\setminus \{1\}$ such that
  $$\max_{s\geq 0}\mathcal E_G(sm_{j_0})<\frac{1}{n}\left(\frac{\alpha_n}{\alpha_0}\right)^{n-1},$$
  where $\alpha_0>0$ is from hypothesis $(f_1)$.
\end{lemma}

{\it Proof.}  By contradiction, we assume that for every $j\in
\mathbb N\setminus \{1\}$, we have
$$\max_{s\geq 0}\mathcal
E_G(sm_{j})\geq\frac{1}{n}\left(\frac{\alpha_n}{\alpha_0}\right)^{n-1}.$$
Since $\mathcal E_G(0)=0$ and $\mathcal E_G(sm_j)\to -\infty$ as
$s\to \infty$ (cf. Lemma \ref{mountain-pass-lemma}), there exists
$s_j>0$ such that $$\max_{s\geq 0}\mathcal E_G(sm_{j})=\mathcal
E_G(s_jm_{j})=s_j^n\frac{\| m_j\|_{0,1}^n}{n}-\mathcal F(s_jm_j).$$

On one hand, since $\mathcal F\geq 0$, the above relations yield
\begin{equation}\label{isten-csuda-1}
    s_j^n\| m_j\|_{0,1}^n\geq
    \left(\frac{\alpha_n}{\alpha_0}\right)^{n-1}.
\end{equation}
 Due to (\ref{m_j-aszimp}),
the above inequality implies that
\begin{equation}\label{utolso-1}
\liminf_{j\to \infty}s_j^n\geq
    \left(\frac{\alpha_n}{\alpha_0}\right)^{n-1}.
\end{equation}

On the other hand,  $s_j>0$ being an extremal point of $s\mapsto
\mathcal E_G(sm_{j}),$ we also have that $\left.\frac{\rm d}{{\rm
d}s}\mathcal E_G(sm_{j})\right|_{s=s_j}=0,$ which is equivalent to
\begin{equation}\label{utolso-2}
s_j^n\| m_j\|_{0,1}^n =\int_Mf(s_jm_j)s_jm_j{\rm d}v_g\ {\rm for\
every}\ j\in \mathbb N\setminus\{1\}.
\end{equation}

By (\ref{m_j-aszimp}), there exists $c_9>0$ such that for large
$j\in \mathbb N$,
\begin{equation}\label{m_j-aszimp-2}
    \|m_j\|_{0,1}^n\leq 1+\frac{c_9}{\log j}.
\end{equation}
Fix
$$L_0>\left(\frac{\alpha_n}{\alpha_0}\right)^{n-1}\omega_n^{-1}e^{c_9\frac{n}{n-1}}.$$
By hypothesis $(f_1),$ there exists $R_1>0$ such that
\begin{equation}\label{f1-hypo-modositas}
    sf(s)e^{-\alpha_0s^\frac{n}{n-1}}\geq L_0 \ \ {\rm for\ every}\
s\geq R_1.
\end{equation}

Note that the sequence  $\{s_j\}$ is bounded. Indeed, if we assume,
up to a subsequence, that   $\lim_{j\to \infty}s_j=\infty$, then for
$j\in \mathbb N$ large enough, we have by (\ref{utolso-2}) that
\begin{eqnarray*}
  \| m_j\|_{0,1}^n
   &\geq& s_j^{-n}\int_{B_{x_0}(\frac{1}{j})}f(s_jm_j)s_jm_j{\rm d}v_g\ \ \ \ \ \ \ \ \ \ \ \ \ \ \ \ \ \ \ \ \ \ \ \ \ \ \ \ \   [sf(s)\geq 0\ {\rm for\ every}\ s\geq 0] \\
   & \geq&
   L_0 s_j^{-n}\int_{B_{x_0}(\frac{1}{j})}e^{\alpha_0(s_jm_j)^\frac{n}{n-1}}{\rm d}v_g\ \ \ \ \ \ \ \ \ \ \ \ \ \ \ \ \ \ \ \ \ \ \ \ \ \ \ \ [{\rm see}\ (\ref{f1-hypo-modositas})]\\
   &=&L_0 s_j^{-n}e^{\alpha_0s_j^\frac{n}{n-1}\omega_{n-1}^{-\frac{1}{n-1}}\log
   j}{\rm
   Vol}_g\left(B_{x_0}\left(\frac{1}{j}\right)\right)\ \ \ \ \ \ \ \ \ \ \ \ \   [{\rm see}\ (\ref{Moser-function-2})]\\
   &\geq&L_0\omega_ne^{n\left(\frac{\alpha_0}{\alpha_n}s_j^\frac{n}{n-1}-1\right)\log
   j-n\log s_j}.\ \ \ \ \ \ \ \ \ \ \ \ \ \ \ \ \ \ \ \ \ \ \ \ \ \ \   [{\rm see\ Proposition\ \ref{comparison}\ (ii)}]
\end{eqnarray*}
Letting $j\to \infty$, on account of (\ref{m_j-aszimp}) we arrive to
a contradiction; thus, $\{s_j\}$ is bounded.



We claim that
\begin{equation}\label{utolso-3}
\lim_{j\to \infty}s_j^n =
    \left(\frac{\alpha_n}{\alpha_0}\right)^{n-1}.
\end{equation}
By contradiction, due to (\ref{utolso-1}), we assume that there
exists $\varepsilon_0>0$ such that (up to a subsequence) for enough
large $j\in \mathbb N$,
$$s_j^\frac{n}{n-1}>\frac{\alpha_n}{\alpha_0}+\varepsilon_0.$$
Note that for every $x\in
    B_{x_0}(\frac{1}{j})$, we have $s_jm_j(x)=s_j{(\log
    j)^\frac{n-1}{n}}{\omega_{n-1}^{-1/n}}\to \infty$  as $j\to \infty$. Therefore, for enough large $j\in \mathbb N,$ relation
    (\ref{f1-hypo-modositas}) can be applied for $s=s_jm_j(x)$ with $x\in B_{x_0}(\frac{1}{j}),$ obtaining in a
    similar manner as above that
$$
s_j^n\| m_j\|_{0,1}^n\geq L_0\omega_ne^{n\left(\frac{\alpha_0}{\alpha_n}s_j^\frac{n}{n-1}-1\right)\log
   j }.
$$
   Consequently, the latter two inequalities, the boundedness of $\{s_j\}$ and (\ref{m_j-aszimp})
   provide a contradiction once $j\to \infty,$ which proves the
   validity of (\ref{utolso-3}).

 For every $j\in \mathbb N\setminus \{1\}$, let $$A_j=\{x\in \overline B_{x_0}(1):s_jm_j(x)\geq R_1\}\ \ {\rm and}\ \ B_j=\overline B_{x_0}(1)\setminus A_j.$$
Due to (\ref{f1-hypo-modositas}), we have
\begin{eqnarray}\label{f-splitting}
  \int_Mf(s_jm_j)s_jm_j{\rm d}v_g &=&\int_{A_j}f(s_jm_j)s_jm_j{\rm d}v_g+\int_{B_j}f(s_jm_j)s_jm_j{\rm
d}v_g  \nonumber\\
   &\geq & L_0 \int_{A_j}e^{\alpha_0(s_jm_j)^\frac{n}{n-1}}{\rm
   d}v_g+\int_{B_j}f(s_jm_j)s_jm_j{\rm
d}v_g\nonumber\\&=&L_0 \int_{\overline
B_{x_0}(1)}e^{\alpha_0(s_jm_j)^\frac{n}{n-1}}{\rm
   d}v_g-L_0 \int_{B_j}e^{\alpha_0(s_jm_j)^\frac{n}{n-1}}{\rm
   d}v_g\nonumber\\&&+\int_{B_j}f(s_jm_j)s_jm_j{\rm
d}v_g.
\end{eqnarray}
Note that $s_jm_j\leq R_1$ in $B_j$, while  $m_j\to 0$ and
$\chi_{B_j}\to 1$ almost everywhere in $\overline B_{x_0}(1)$ as
$j\to \infty.$ Consequently, on one hand, by the Lebesgue dominated
convergence theorem we have
$$\lim_{j\to \infty}\int_{B_j}e^{\alpha_0(s_jm_j)^\frac{n}{n-1}}{\rm
   d}v_g=\int_{\overline B_{x_0}(1)}{\rm
   d}v_g={\rm Vol}_g(\overline B_{x_0}(1))\ \ {\rm and} \ \ \lim_{j\to \infty}\int_{B_j}f(s_jm_j)s_jm_j{\rm
d}v_g=0.$$ On the other hand,
$$\int_{\overline
B_{x_0}(1)}e^{\alpha_0(s_jm_j)^\frac{n}{n-1}}{\rm
   d}v_g=\int_{\overline
B_{x_0}(1)\setminus
B_{x_0}(\frac{1}{j})}e^{\alpha_0(s_jm_j)^\frac{n}{n-1}}{\rm
   d}v_g+\int_{B_{x_0}(\frac{1}{j})}e^{\alpha_0(s_jm_j)^\frac{n}{n-1}}{\rm
   d}v_g=:I_j^1+I_j^2.$$
Clearly, we have $I_j^1\geq 0$, and for large $j\in \mathbb N,$
\begin{eqnarray*}
  I_j^2 &=& \int_{B_{x_0}(\frac{1}{j})}e^{\alpha_0(s_jm_j)^\frac{n}{n-1}}{\rm
   d}v_g \\
   &\geq& \int_{B_{x_0}(\frac{1}{j})}e^{\alpha_nm_j^\frac{n}{n-1}\|m_j\|_{0,1}^{-\frac{n}{n-1}}}{\rm
   d}v_g=e^{n\log j\|m_j\|_{0,1}^{-\frac{n}{n-1}}}{\rm Vol}_g\left(B_{x_0}(\frac{1}{j})\right) \ \ \ \ \ \ \ \ \ \ \ \ \ \ \    [{\rm see}\ (\ref{isten-csuda-1})]\\
   &\geq&\omega_n j^{n\left(\|m_j\|_{0,1}^{-\frac{n}{n-1}}-1\right)}\ \ \ \ \ \ \ \ \ \ \ \ \ \ \ \ \ \ \ \ \ \ \ \ \ \ \ \ \ \ \ \ \ \ \ \ \ \ \ \ \ \ \ \ \ \ \ \ \ \ \  \ \   [{\rm see\ Proposition\ \ref{comparison}\
   (ii)}]\\
   &\geq& \omega_n j^{n\left(\left(1+\frac{c_9}{\log j}\right)^{-\frac{1}{n-1}}-1\right)}.\ \ \ \ \ \ \ \ \ \ \ \ \ \ \ \ \ \ \ \ \ \ \ \ \ \ \ \ \ \ \ \ \ \ \ \ \ \ \ \ \ \ \ \ \  \ \ \ \ \ \ \ \ \ \ \ \ \ \ \ \ \ \  \ \
     [{\rm see\  (\ref{m_j-aszimp-2}})]
\end{eqnarray*}
Therefore, $$\liminf_{j\to \infty}I_j^2\geq \omega_n\lim_{j\to
\infty} j^{n\left(\left(1+\frac{c_9}{\log
j}\right)^{-\frac{1}{n-1}}-1\right)}=\omega_ne^{-c_9\frac{n}{n-1}}.$$
Putting in (\ref{utolso-2}) the latter estimates together with
relations (\ref{m_j-aszimp}), (\ref{utolso-3}) and
(\ref{f-splitting}), it follows that
$$\left(\frac{\alpha_n}{\alpha_0}\right)^{n-1}\geq L_0\omega_ne^{-c_9\frac{n}{n-1}},$$
which contradicts the choice of $L_0.$ The proof is complete. \hfill
$\square$\\

{\it Proof of Theorem \ref{theorem-existence}.} Let $m_{j_0}\in
W_G^{1,n}(M)$ be the  Moser function which satisfies the conclusion
of Lemma \ref{Moser-fuggvenyek-fo}. By Lemma
\ref{mountain-pass-lemma}, the functional $\mathcal
E_G:W_G^{1,n}(M)\to \mathbb R$ has the mountain pass geometry; in
particular, if $e_0=s_0 m_{j_0}\in W_G^{1,n}(M)$ with $s_0>0$ large
enough, then $\mathcal E_G(e_0)<0=\mathcal E_G(0)$ and $\mathcal
E_G(u)\geq \tilde \delta>0$ for every $u\in W_G^{1,n}(M)$ with
$\|u\|_{0,1}=\tilde r$, where $\tilde r<\|e_0\|.$ By using the
mountain pass lemma for $\mathcal E_G$ without the Palais-Smale
compactness condition, see e.g. Brezis and Nirenberg \cite[p.
943]{BN}, there exists a sequence $\{u_j\}\subset W_G^{1,n}(M)$ such
that
\begin{equation}\label{energia-vege}
    \mathcal E_G(u_j)\to c\ \ {\rm and}\ \ \mathcal E_G^\prime(u_j)\to
0\ {\rm in}\ W_G^{1,n}(M)^*,
\end{equation}
where $$c=\inf_{\lambda\in \Lambda}\max_{s\in [0,1]}\mathcal
E_G(\lambda(s))\geq \tilde \delta,$$ and $\Lambda=\{\lambda\in
C([0,1],W_G^{1,n}(M)): \lambda(0)=0,\ \lambda(1)=e_0\}.$ According
to Lemma \ref{Palais-Smale}, there exists $u_G\in W_G^{1,n}(M)$ such
that, up to a subsequence,
\begin{equation}\label{FuG-konvergencia}
    \lim_{j\to \infty}\mathcal F(u_j)= \mathcal F(u_G),
\end{equation}
$u_j\to u_G$ strongly in $L^p(M)$ for every $p\in (n,\infty),$ and
$u_G$ is a
  critical point of $\mathcal E_G$. The latter fact  with
  Lemma \ref{lemma-szimmetria-elv} shows that $u_G$ is a non-negative $G-$invariant weak solution of $(\mathcal P)$.

It remains to prove that $u_G\neq 0.$ By contradiction, if $u_G=0$,
relations (\ref{energia-vege}) and (\ref{FuG-konvergencia}) imply on
one hand that
\begin{equation}\label{ujG}
\lim_{j\to \infty}\|u_j\|_{0,1}^n=nc\geq n\tilde \delta>0.
\end{equation}
On the other hand, if we apply $u_j$ as a test-function in $\mathcal
E_G^\prime(u_j)\to 0$, one has $\lim _{j\to \infty}\langle\mathcal
E_G^\prime(u_j),u_j\rangle_*=0$, i.e.,
\begin{equation}\label{micsolda-marha}
    \lim_{j\to \infty}\left(\|u_j\|_{0,1}^n-\int_M f(u_j)u_j{\rm
d}v_g\right)=0 .
\end{equation}
 By the definition of the minimax value $c$, it
follows by Lemma \ref{Moser-fuggvenyek-fo} that
$$c\leq \max_{s\in [0,1]}\mathcal
E_G(se_0)\leq \max_{s\geq 0}\mathcal
E_G(sm_{j_0})<\frac{1}{n}\left(\frac{\alpha_n}{\alpha_0}\right)^{n-1}.$$
This estimate and  (\ref{ujG}) guarantee the existence of $q>n$ such
that for every large $j\in \mathbb N,$
$$\frac{q}{q-1}\|u_j\|_{0,1}^\frac{n}{n-1}<\frac{\alpha_n}{\alpha_0}.$$
Relations (\ref{f-elso-becsles}), (\ref{Phi-egyenlotlenseg}), the
 H\"older's inequality and the latter relation imply that for large
$j\in \mathbb N$,

   \begin{eqnarray*}
  0\leq \int_M f(u_j)u_j{\rm d}v_g &\leq& c_0\int_M  |u_j|^{\gamma}{\rm d}v_g +c_0\int_M|u_j|\Phi_n(\alpha_0|u_j|^\frac{n}{n-1}){\rm d}v_g \\
     &\leq& c_0\|u_j\|_{L^\gamma(M)}^\gamma+c_0\|u_j\|_{L^q(M)}\left(\int_M\Phi_n\left(\alpha_0\frac{q}{q-1}|u_j|^\frac{n}{n-1}\right){\rm d}v_g \right)^{1-\frac{1}{q}}\\
     &\leq& c_0\|u_j\|_{L^\gamma(M)}^\gamma+c_0\|u_j\|_{L^q(M)}\left(S_{\alpha_n,1}^0(M,g)
     \right)^{1-\frac{1}{q}}.
   \end{eqnarray*}
Note that  $S_{\alpha_n,1}^0(M,g)<\infty$, cf. Yang, Su and Kong
\cite[Theorem 1.2]{YSK}. Moreover, since $\gamma,q>n$ and
$\lim_{j\to \infty}\|u_j\|_{L^p(M)}= 0$ for every $p\in (n,\infty)$,
it follows from the last estimate that
$$\lim_{j\to \infty}\int_M f(u_j)u_j{\rm d}v_g=0.$$
This limit and relations (\ref{ujG}) and (\ref{micsolda-marha})
provide a contradiction. Therefore, $u_G\neq 0$. \hfill $\square$

\section{Examples and an open problem}

 \noindent We present some possible scenarios where
Theorem \ref{theorem-existence} can be applied.

\begin{example}\label{example-1} {\rm [{\it Euclidean case}]} \rm  If $(M,g)=(\mathbb R^n,g_{\rm euc})$ is the
usual Euclidean space, Theorem \ref{theorem-existence} can be
applied  for $x_0=0$ and $G={\sf SO}(n_1,\mathbb R)\times...\times
{\sf SO}(n_l,\mathbb R)$ with $n_j\geq 2$,  $j=1,...,l$ and
$n_1+...+n_l=n$, where ${\sf SO}(m,\mathbb R)$ is the special
orthogonal group in $\mathbb R^m$. Indeed, we have ${\rm
Fix}_{\mathbb R^n}(G)=\{0\}$ and $O_G^x=|x_{n_1}|\mathbb
S^{n_1-1}\times...\times |x_{n_l}|\mathbb S^{n_l-1}$ for each
$x=(x_{n_1},...,x_{n_l})\in \mathbb R^{n_1}\times...\times \mathbb
R^{n_l}\setminus\{0\}.$ \hfill $\square$
\end{example}

\begin{example} {\rm [{\it Hyperbolic case}]} \rm  For the hyperbolic space we use the
Poincar\'e ball model $\mathbb H^n=\{x\in \mathbb R^n:|x|<1\}$
endowed with the Riemannian metric $g_{\rm
hyp}(x)=(g_{ij}(x))_{i,j={1,...,n}}=\frac{4}{(1-|x|^2)^2}\delta_{ij}$.
It is well known  that $(\mathbb H^n,g_{\rm hyp})$ is a homogeneous
Hadamard manifold with constant sectional curvature $-1$. Theorem
\ref{theorem-existence} can be applied with the same choice for
$x_0$ and $G$ as in Example \ref{example-1}.\hfill $\square$
\end{example}

\begin{example} {\rm [{\it Symmetric positive definite matrices}]}
\rm Let ${\rm Sym}(n,\mathbb R)$ be the set of symmetric $n\times n$
matrices with real values, ${\rm P}(n,\mathbb R)\subset {\rm
Sym}(n,\mathbb R)$ be the $\frac{n(n+1)}{2}-$dimensional cone of
symmetric positive definite matrices, and ${\rm P}(n,\mathbb R)_1$
be the subspace of matrices in ${\rm P}(n,\mathbb R)$ with
determinant one.
 The set ${\rm P}(n,\mathbb R)$ is endowed with the scalar product
$$
  \langle\langle U,V \rangle\rangle_X={\rm Tr}(X^{-1}VX^{-1}U)\ \ {\rm for\ all}\ \ X\in
{\rm P}(n,\mathbb R),\ U,V\in T_X({\rm P}(n,\mathbb R))\simeq {\rm Sym}(n,\mathbb R), $$
where ${\rm Tr}(Y)$ denotes the trace of $Y\in {\rm Sym}(n,\mathbb
R)$, and let us denote by $d_H:{\rm P}(n,\mathbb R)\times {\rm
P}(n,\mathbb R)\to \mathbb R$ the induced metric function.  The pair
$({\rm P}(n,\mathbb R),\langle\langle\cdot,\cdot \rangle\rangle)$ is
a Hadamard manifold, see Lang \cite[Chapter XII]{Lang}. Note that
${\rm P}(n,\mathbb R)_1$ is a convex totally geodesic submanifold of
${\rm P}(n,\mathbb R)$ and the special linear group  ${\sf
SL}(n,\mathbb R)$ leaves ${\rm P}(n,\mathbb R)_1$ invariant and acts
transitively on it; thus $({\rm P}(n,\mathbb
R)_1,\langle\langle\cdot,\cdot \rangle\rangle)$ is itself a
homogeneous Hadamard manifold, see Bridson and Haefliger
\cite[Chapter II.10]{BH}. Moreover, for every $\sigma\in {\sf
SL}(n,\mathbb R)$, the map
 $[\sigma]:{\rm P}(n,\mathbb
R)_1\to {\rm P}(n,\mathbb R)_1$ defined   by $[\sigma](X)= \sigma
X\sigma^t$, is an isometry; here, $\sigma^t$ denotes the
 transpose of $\sigma.$ 

  Let $G={\sf SO}(n,\mathbb R)$. One can prove that ${\rm Fix}_{{\rm
P}(n,\mathbb R)_1}(G)=\{I_n\}$, where $I_n$ is the identity
 matrix. First, it is clear that $I_n\in {\rm Fix}_{{\rm
P}(n,\mathbb R)_1}(G)$; indeed, for every $\sigma\in G$ we have
$[\sigma](I_n)=\sigma I_n\sigma^t=\sigma \sigma^t=I_n.$ Second, if
$X_0\in {\rm Fix}_{{\rm P}(n,\mathbb R)_1}(G)$, then it turns out
that $\sigma X_0=X_0\sigma$ for every $\sigma\in G$. By using
elementary matrices from $G$, the latter relation implies that
$X_0=c I_n$ for some $c\in \mathbb
 R$. Since $X_0\in {\rm P}(n,\mathbb R)_1$, we necessarily have
 $c=1.$ Moreover, the orbit $O_G^X$ of the matrix $X\in {\rm P}(n,\mathbb
R)_1\setminus \{I_n\}$ under the action of $G$ is the geodesic
sphere in ${\rm P}(n,\mathbb R)_1$ with center $I_n$ and radius
$d_H(I_n,X);$ in particular, Card$(O_G^X)=\infty.$ Indeed, for every
$\sigma\in G$, since $[\sigma]$ is an isometry on ${\rm P}(n,\mathbb
R)_1$, it follows that
\begin{eqnarray*}
 d_H^2(I_n,[\sigma](X)) &=& d_H^2([\sigma](I_n),[\sigma](X))=  d_H^2(I_n,X).
\end{eqnarray*}
Consequently, Theorem \ref{theorem-existence} is applicable on ${\rm
P}(n,\mathbb R)_1$ with the choices $x_0=I_n$ and $G={\sf
SO}(n,\mathbb R)$, respectively.
\end{example}

We conclude the paper with the following open problem concerning the
 volume growth of geodesic balls in the presence of the  Moser-Trudinger
inequality $({\bf MT})_{\alpha,1}^0$:\\

\noindent  {\bf Problem.} {\it Let $(M,g)$ be an $n-$di\-men\-sional
complete non-compact Riemannian manifold $(n\geq 2)$ with
non-negative Ricci curvature and assume the Moser-Trudinger
inequality $({\bf MT})_{\alpha,1}^0$ holds on $(M,g)$ for some
$\alpha\in (0,\alpha_n]$. Is there any $\gamma>0$ such that $${\rm
Vol}_g(B_x(r))\geq
\left(\frac{\alpha}{\alpha_n}\right)^\gamma\omega_n r^n\ \ {for\
every}\ x\in M\ and \ r>0?$$}

If the answer is affirmative, we could state that the sharp
Moser-Trudinger inequality $({\bf MT})_{\alpha_n,1}^0$ holds on  an
$n-$di\-men\-sional complete non-compact Riemannian manifold $(M,g)$
with non-negative Ricci curvature if and only if $(M,g)$ is
isometric to the Euclidean space $\mathbb R^n$. Similar results can
be found e.g. in do Carmo and Xia \cite{doCarmo-Xia}, Krist\'aly \cite{K-Potential},  Krist\'aly and
Ohta \cite{Kri-Ohta}, Ledoux \cite{Ledoux-CAG} and references
therein for various Sobolev-type inequalities; the arguments in
these papers are based on the precise shape of extremal functions
for the studied Sobolev-type inequalities in the Euclidean setting.
Although Li and Ruf \cite{LR} proved that the supremum $S_n^{LR}$ in
(\ref{Li-Ruf-eredmeny}) is achieved, no explicit extremal function
is known.\\

 \noindent {\bf Acknowledgments.} The author thanks Professor
Philippe G. Ciarlet for his invitation to the City University of
Hong Kong where this work has been initiated,  Professor Yunyan
Yang for stimulating discussions on his papers \cite{Adimurthi-Yang,
doO-Yang, Yang-JFA2012}, and the anonymous 
Referee
for
her/his valuable
comments.

\end{document}